\documentclass[authoryear,preprint,review,12pt]{elsarticle}

\usepackage{apalike} % reference format

\usepackage{amssymb}
\usepackage{amsmath}

\newtheorem{thm}{Theorem}
\newtheorem{lem}{Lemma}

\newtheorem{cor}{Corollary}

\newenvironment{proof}{{\it Proof:\quad}}{\hfill $\blacksquare$\par}
\newcommand{\argmin}{\mathop{\mathrm{argmin}}}
\newcommand{\argmax}{\mathop{\mathrm{argmax}}}

\journal{arXiv.org}
\begin{document}

\begin{frontmatter}



%%%%%%%%%%%%%%%%%%%%%%
% Synonyms: relieve, alleviate, lighten, assuage, mitigate, allay mean to make something less grievous.

% relieve implies a lifting of enough of a burden to make it tolerable.
% took an aspirin to relieve the pain

% alleviate implies temporary or partial lessening of pain or distress.
% the lotion alleviated the itching

% mitigate suggests a moderating or countering of the effect of something violent or painful.
% the need to mitigate barbaric laws 

%\title{Fine Control of Conservatism for Robust Optimization by Adjustable Regret}
\title{Mitigate Overconservatism for Robust Optimization by Adapting to Opportunities}


\author[inst1]{Yingjie Lan (email: ylan@pku.edu.cn)}

\affiliation[inst1]{organization={Peking University HSBC Business School},%Department and Organization
            addressline={Xili University Town}, 
            city={Shenzhen},
            postcode={518005}, 
            state={Guangdong},
            country={China}}

%\author[inst2]{Author Two}
%\author[inst1,inst2]{Author Three}

%\affiliation[inst2]{organization={Department Two},%Department and Organization
%            addressline={Address Two}, 
%            city={City Two},
%            postcode={22222}, 
%            state={State Two},
%            country={Country Two}}

\begin{abstract}
Overconservatism has long been recognized as a major issue with robust optimization, despite its key advantages of tractability, performance guarantee, and limited information.
To address this issue, a new criterion is proposed that can adapt its level of conservatism continuously to the opportunities out there, while maintaining all the key advantages just mentioned. %gamut
With this criterion, a general framework of  conservatism control based on optimal performance guarantee is developed and characterized, and a new approach to competitive ratio analysis is established.
The criterion is then applied to the robust one-way trading problem, where analytical solution is obtained, and the competitive ratio is derived directly via the new approach.
Numerical experiments are conducted to demonstrate the effectiveness of conservatism control based on the new criterion, with the average reward improvable by 4\% - 17\% over the other commonly used criteria. 
\end{abstract}

%without tampering with the uncertainty set.
%and inclusion of various commonly used criteria at special values of the control parameter.
%Tractability for robust linear programs with this criterion is established by reformulating them into those with the maximin criterion, for which tractable solution schemes and theoretical results are actively developed in the literature.




\end{frontmatter}

\section{Introduction.}
Robust optimization (RO) is a popular method of decision-making under uncertainty, where the decision criterion plays a key role in gaining the main advantages of RO: limited information, performance guarantee, and computational tractability.
Limited information refers to the uncertainty set of all scenarios without probabilistic distributions as needed in stochastic optimization. 
% often represented as realizations of parameters,
% The uncertainty set can be inferred from expert judgments or limited data (\citealt{goerigk2023data}).
As the outcome of an action also depends on the realized scenario, it is impossible to evaluate an action under limited information without a proper criterion.
%, As the outcome of an action depends on the realized scenario. 
% Without probabilities over the uncertainty set, it is impossible to compute expectations for the objective that is used in stochastic optimizaiton, thus a proper decision criterion must be employed to form an objective that evaluates a solution under all scenarios, which is key to model with limited information. 
% The criterion also plays a crucial role for performance guarantee and computational tractability.
Performance guarantee comes from the criterion considering the worst scenarios when evaluating actions, which brings about performance robustness.
% the recommended solutions are robust as
And finally, computational tractability rests on the criterion preserving convexity, as many classes of convex optimization problems admit polynomial-time algorithms (\citealt{nemirovsky1994interior}), while mathematical optimization is NP-hard in general (\citealt{murty1987some}).
Since convexity of the feasible set is preserved under set intersection, the RO model remains convex as long as the criterion produces a convex objective. % as a result of intersection of convex sets as feasible region for each scenario
% It is even better if linearity, the simplest form of convexity, can be preserved.
Criteria that uphold these advantages, such as maximin profit, minimax cost, or minimax regret, have made RO an appealing method in numerous fields of applications including finance, operations management, aerospace, robotics, and defense.

Despite such advantages, a major issue with RO is overconservatism, %which has led to flurries of fruitful researches.
which sacrifices too much performance for the sake of robustness, due to an obsession with the worst scenarios while ignoring all other opportunities, no matter how likely they may occur.
This performance issue seriously hinders the adoption of RO in some industries, such as robust revenue management in airlines (see \citealt{vinod2021approach}), where the profit margin is already razor-thin so that the performance sacrificed is often a dealbreaker.
When such obsession is assuaged by having some ambiguous distributional information on uncertainty, the distributionally robust optimization (DRO) makes a less conservative methodology, see \cite{kuhn2025distributionally} for a recent survey and \cite{zhen2025unified} for a unified theory.
Meanwhile, tackling overconservatism without distributions has always been an important and active research front ever since the first RO models appear in \cite{soyster1973convex}. 
%for linear programming (LP).
% It often comes by extreme pessimism, which expects things to go wrong with negative outcomes, leading to a risk-averse attitude to protect against worst outcomes.
The earliest models are extremely pessimistic on how uncertainty affects feasibility and performance, which leads to overconservatism in two distinct ways: 
Extreme pessimism on feasibility insists that a solution must be feasible for all scenarios, which may end up excluding some good solutions as infeasible; 
% though considering all scenarios delivers solid feasibility robustness
In a different way, extreme pessimism on performance demands that a solution shall be evaluated by its performance in the worst scenario, while ignoring all other scenarios that may offer great opportunities.
Both ways can take a toll on performance, thus overconservatism has been tackled in the literature by either excluding extreme scenarios from feasibility test, or resorting to less conservative criteria for solution evaluation.

%as uncertainty can cause violations of the constraints. 

%The guarantees in RO that offer robustness comes at the sacrifice of performance, which leads to overconservatism if too much is sacrificed.
%While RO can provide robust solutions to real-world problems, it also leads to over-conservative solutions that deliver relatively poor performances. 

%The key lies in balancing robustness with optimality, which is one of the ongoing challenges in this field.

Scenario exclusion for conservatism reduction begins by questioning the absolute feasibility guarantee under all scenarios. 
Though it is necessary for critical applications where infeasibility causes disasters like doomed satellites or ruined rovers, such absolute guarantee can be relaxed if adverse events only bring about limited consequences, such as low demand or supply in business applications.
In the latter case, it is acceptable to exclude some rare and extreme scenarios from the uncertainty set, settling for a probabilistic guarantee in exchange for better performances.
Researches in this regard provide insights into robust solutions as well as probabilistic guarantees of feasibility, see \cite{ben2008selected} and \cite{bertsimas2011theory} for a comprehensive survey and \cite{ben2009robust} for a book treatment. 
%enrich theories of RO (\citealt{ben2009robust,bertsimas2011theory,gabrel2014recent}),
Models with ellipsoidal uncertainty sets are proposed by %Ben-Tal and Nemirovski 
\cite{ben1998robust,ben1999robust,ben2000robust},  %El-Ghaoui et al. 
\cite{el1997robust}, and \cite{el1998robust}. 
As such models are nonlinear and computationally demanding, \cite{bertsimas2004price} proposes uncertainty budget to control the uncertainty set while maintaining linearity.
Scenario exclusion may pose practical challenges in estimating the cost of infeasibility and probability for those rare and extreme scenarios to make balanced trade-offs between robustness and performance.

%%Brief review on different criteria for RO
Alternative criteria have also been developed as another way to tackle overconservatism by the researchers.
The minimax criterion adopted in the earliest RO models is invented when economists contemplated on decision theories (\citealt{wald1950statistical}).
This criterion is criticized by \cite{savage1954foundations} as ``ultrapessimistic", as it focuses entirely on performances in the worst scenario and ignores all plausible opportunities in others. 
By incorporating opportunity loss as regret, \cite{savage1951theory} proposes a less conservative criterion that minimizes the worst-case regret, known as the minimax regret criterion that provides a {\em regret guarantee}.
After an action is taken and then a scenario is realized, the decision maker may regret not taking the ex post optimal action to fully exploit the opportunities offered by the realized scenario. 
The regret is the opportunity loss as measured by the difference between the ex post optimal and the actual objective.
% The opportunity in a scenario is fully realized in the ex post optimal objective, which may serve as a measure of the opportunity. 
%where the regret is the performance loss compared to the ex post optimal.
%it is great to leave as little untaken as possible.   
The so-called ``competitive ratio," a popular criterion for online optimization (\citealt{borodin2005online, kouvelis2013robust}), is equivalent to the relative regret criterion, which considers the ratio of regret to the ex post optimal objective instead for a relative performance guarantee.
Various numerical studies find the absolute regret less conservative than the relative regret, which is in turn less conservative than the maximin criterion in reward maximizing problems (\citealt{lan2008revenue,perakis2008regret, poursoltani2021adjustable}). %for maximization problems
% (\citealt{
% perakis2008regret,
% natarajan2014probabilistic, %Natarajan et al. 2014, 
% wang2016competitive,
% caldentey2017intertemporal, %Caldentey et al. 2017). 
% poursoltani2021adjustable}),

%\cite{poursoltani2021adjustable} even consider it as the Achilles' heel of traditional RO, for it focuses entirely on worst-case profits and ignores all plausible opportunities for higher profits, leading to highly conservative solutions. 

These commonly used criteria maintain the key advantages of RO and offer three distinct levels of conservatism, but for fine control of conservatism, it would require a continuum of choices.
% present the reasons for fine grained control of conservatism, and the right match of conservatism with the circumstances of application
An early attempt is made by \cite{hurwicz1951optimality} with the Hurwicz criterion, 
which evaluates a solution by a weighted sum of its worst and best outcomes in order to balance pessimism and optimism. 
There is a ``coefficient of pessimism" ($\alpha$) between $0$ and $1$ as the weight on the worst outcome, while ($1-\alpha$) weights on the best outcome.
% For cost minimization, 
It becomes the most conservative as the minimax criterion when $\alpha=1$, and the most aggressive as the minimin criterion when $\alpha=0$, % neither extremes are good
and in between it generally gets more and more conservative as $\alpha$ increases. 
Unfortunately, it can not preserve convexity and maintain computational tractability. 
%and usually poses great challenges for theoretical analysis.
%Analysis by \cite{gaspars2014modifications} also shows that sometimes it can lead to quite illogical answers, thus its behavior can be difficult to understand sometimes.

% This work aspires to create a new criterion for fine control of conservatism while keeping all the major advantages of RO. %intact.
Many other criteria are proposed after Hurwicz for fine control of conservatism, but none of them can keep the key advantages of RO. %overcome the disadvantage of lost convexity. to moderate conservatism 
For example, the $p$-robustness by \cite{snyder2006facility} first screens out by additional constraints overly conservative solutions whose worst-case regret exceeds an upper limit, but sometimes it is very difficult to determine if it is making a problem infeasible.
%\cite{roy2010robustness} proposes the $bw$-robustness criterion, which guarantees an objective value of at least $w$ in all scenarios, and maximizes the probability of reaching a target value of $b (b > w)$ in the presence of a known distribution of parameters. 
% Kalai, Lamboray, and Vanderpooten
\cite{kalai2012lexicographic} suggest the lexicographic robustness criterion to mitigate the primary role of the worst-case scenario in solution evaluation, yet it requires finite uncertainty sets, and poses a serious computational challenge for large uncertainty sets.

To fill this gap, adjustable regret minimization (ARM) is proposed as a new criterion with a continuum of conservatism choices for fine control of conservatism, while upholding all the key advantages of RO.
The ARM criterion provides guarantees under limited information by construction, and maintains tractability by preserving convexity.
A general framework for conservatism control by the ARM criterion is developed to improve performance by adapting to opportunities.
%as demonstrated by applying it to the one-way trading problem with numerical studies.

%TODO: motivate most likely opportunities here: 1. there are applications for this. 2. a little extra information can help reduce the obsession on worst scenarios.
The framework of conservatism control follows a common sense: Blind pessimism or optimism is ineffective or even detrimental to performance as they fail to account for reality, and high performance depends on a realistic attitude that adapts to the opportunities out there. 
It makes use of likely opportunities or typical scenarios, which may be provided by experts, such as the nominal values used to construct the uncertainty set.
A nominal value is often the point forecast --- the single best estimate provided by data or domain experts, corresponding to the mean or mode of an uncertain parameter's distribution, which is a common practice in project management as experts are asked to provide the most likely duration of tasks.
%In practice data often clusters around a single typical value, with values becoming less frequent when deviating further from it, which is known as unimodality.
The level of conservatism is then chosen from the continuum offered by the ARM criterion to catch such opportunities with an optimal performance guarantee, so that some sacrificed performance may be salvaged.
%by which a mechanism is devised to choose from the continuum of conservatism offered by the ARM criterion.
On a side note, the ARM criterion and the framework can be adapted to DRO if the likely distribution in the ambiguity set is known, such as the nominal distribution.

The ARM criterion also supports a new approach to competitive ratio analysis that may significantly reduce the analysis complexity to derive closed-form solutions.
Competitive ratio analysis is often more complex than absolute regret analysis, which can be observed in comparing \cite{el2001optimal} and \cite{wang2016competitive}, as each carries out one type of analysis for exactly the same problem setup.
This new approach first solves the problem with the ARM criterion in similar complexity as the absolute regret analysis, then the competitive ratio is derived directly by solving an equation established by theory, which is applied to robust one-way trading later.
Conceivably, algorithms for numerical computation of competitive ratios can also be designed from this approach in the future. 
% For this purpose, the ARM criterion will be introduced and studied under a general multi-stage setting, which includes the single-stage and two-stage as special cases. 

% that offers smooth control of conservatism by  interpolation with good analytical properties,  gives rise to a new approach to competitive ratio analysis that may reduce analysis complexity.

The contributions of this paper are as follows:
(i) The ARM criterion is proposed for fine control of conservatism, with its properties studied theoretically to confirm that the major advantages of RO are maintained. % without altering the uncertainty set
%The properties of the ARM criterion are studied theoretically, such as convexity preservation and conservatism control.
% It does not tamper with the uncertainty set, and works with any uncertainty set. 
% It can be applied either independently or jointly with methods that modify the uncertainty set. 
% and the ARM criterion is proven to maintain linearity for LP problems with certain structures. 
(ii) A framework based on optimal performance guarantee for conservatism control is devised with its mechanism studied, so that a level of conservatism may be chosen from the continuum offered by the ARM criterion to adapt to likely opportunities.
% Linear RO problems with the ARM criterion can be reformulated as traditional RO problems, and can therefore take advantage of the tractable solution schemes and theoretical results actively developed recently for this class of problems. 
(iii) A new approach to competitive ratio analysis based on the ARM criterion is introduced, which can reduce the complexity  to find analytical solutions, or help design more efficient numerical procedures.
(iv) The robust one-way trading problem with the ARM criterion is solved analytically, and the competitive ratio directly derived by the new approach. Numerical experiments are conducted to demonstrate the effectiveness of conservatism control, with the average reward improvable by 4\% to 17\% over other commonly used criteria. %and that a different problem setting calls for a different control parameter of conservatism in the ARM criterion. 
% Some interesting observations are made: 1. overly conservative solutions could actually suffer higher overall risk and lower average reward at the same time. 2. different problem settings could enjoy different degree of conservatism, which may go beyond interpolation and call for extrapolation.
	
The rest of this paper is organized as follows.
Section \ref{sec:form} gives general ARM formulations, whose properties are studied in Section \ref{sec:prop}, leading to a new approach to competitive ratio analysis. 
The framework and its mechanism for conservatism control is also developed and analyzed.
%Then section \ref{sec:tractability} deals with the tractability of linear problems, 
In Section \ref{sec:1-way} the ARM criterion is applied to the robust one-way trading problem. 
The analysis derives a closed-form solution, which readily yields the competitive ratio by the new approach.
And numerical experiments are conducted to demonstrate effective control of conservatism by the ARM criterion.
Finally, Section \ref{sec:conclusion} draws conclusions with future research outlooks.

\section{Formulations.}\label{sec:form}

For a gentle introduction, the ARM criterion is first formulated with single-stage problems, and extended later to multi-stage problems.
A scenario $\zeta$ consists of realized values for uncertain data, and the uncertainty set $\mathcal{U}$ has all scenarios.
% The set $\mathcal{U}$ is assumed a closed set, as it is so in most practical and theoretical settings. 
Let $X_\zeta$ be the feasible set for $\zeta$, as the constraints may depend on $\zeta$. 
The set $X=\bigcap_{\zeta \in \mathcal{U}} X_\zeta$ is robustly feasible for all scenarios. 
To be general, both $\mathcal{U}$ and $X_\zeta$ may be continuous or discrete.
%The set $X$ may also be continuous or discrete, for example, $X=\{ x\in X: \forall \zeta \in \mathcal{U}, g(x; \zeta)\le 0 \}$, where $g(x;\zeta): X\times\mathcal{U} \rightarrow \mathbb{R}^{n_g}$, and $X=\mathbb{R}^{n_r}\times\mathbb{Z}^{n_z}$ is a mixture of continuous and discrete space with $n_r, n_z \ge 0$ being the number of continuous and discrete components in $x$. 
In single-stage problems, an action $x$ usually takes place first, then a scenario $\zeta$ realizes, and the reward $r(x,\zeta)$ depends on both $x$ and $\zeta$. 
Let $r^*(\zeta) = \max_{x\in X_\zeta} r(x,\zeta)$ denote the ex post optimal, which measures the opportunities in $\zeta$.
It is assumed throughout that the `$\min$' and `$\max$' operators are well-defined, otherwise `$\inf$' and `$\sup$' can be used instead.

The $\beta$-adjusted regret $D(x, \zeta; \beta) = \beta r^*(\zeta) - r(x,\zeta)$ is considered in the ARM criterion for an action $x\in X$ after a scenario $\zeta$ is realized, where $\beta \in [-\infty, +\infty]$ is the {\em conservatism control parameter} (CCP), a constant multiplier in the benchmark $\beta r^*(\zeta)$.
The worst regret $\bar{D}(x;\beta) = \max_{\zeta \in \mathcal{U}} D(x, \zeta; \beta)$ serves as a regret guarantee provided by $x$, which is minimized for a recommended action that provides the best regret guarantee $D(\beta)$:
\begin{equation}\label{defn:arm-beta}
	D(\beta) = \min_{x\in X} \bar{D}(x;\beta)
	= \min_{x\in X} \max_{\zeta\in \mathcal{U}} \beta r^*(\zeta) - r(x, \zeta),
\end{equation}
% \begin{eqnarray}\label{defn:arm-beta}
% D(\beta) & = & \min_{x\in X} \bar{D}(x;\beta)
% = \min_{x\in X} \max_{\zeta\in \mathcal{U}} \beta r^*(\zeta) - r(x, \zeta),\nonumber\\
% &=& \min_{x\in X} \max_{\zeta\in \mathcal{U}} \beta \{\max_{x'\in X} r(x', \zeta)\} - r(x, \zeta).
% \end{eqnarray} 
%Note that $D(x, \zeta; \beta=1)$ gives the traditional definition of absolute regret.
from which the {\em performance guarantee} provided by an optimal solution for a $\zeta$ with $r^*(\zeta)=r$ is simply $\check{r}(r, \beta) = \beta r - D(\beta)$.
%If experts can provide the most likely opportunities in terms of $r^*(\zeta)$, the performance guarantee may be maximized to best catch them, which will be elaborated upon later.
Note that though \eqref{defn:arm-beta} has been used for reducing computational complexities in combinatorial optimization (\citealt{averbakh2005computing}), it has never been proposed as a criterion for moderating conservatism.
Observe that the ARM criterion with a given $\beta$ would recommend the same solutions if the reward $r(x,\zeta)$ is replaced by $r'(x,\zeta)= k r(x,\zeta) + b$ for any $k>0$ and $b\in \mathbb{R}$.
Also note that it is applicable to DRO if $\zeta$ represents a distribution and $r(x,\zeta)$ is the expected reward.

A glimpse of the effects of $\beta$ on conservatism may be taken as the ARM criterion transmorphs into other well-known criteria with different $\beta$ values. % as $\beta$ takes on values in $[0,\infty)$. %  that are unified into a continuum
At $\beta=0$ it becomes the maximin criterion in traditional RO models, which is the most conservative. 
When $\beta$ takes on the competitive ratio (usually a special value between 0 and 1), it is the same as the relative regret criterion (more details later), which is less conservative than the previous.
At $\beta = 1$ it is the absolute regret criterion that is even less conservative, and finally as $\beta = +\infty$ it recovers the maximax criterion that is the most aggressive.
Interestingly, at the other extreme of $\beta=-\infty$ is the opposite to the maximax criterion, which recommends an $x$ that performs best for the least promising scenario $\zeta'$ with $r^*(\zeta')=\min_{\zeta\in U} r^*(\zeta)$.
These special cases suggest that the ARM criterion becomes increasingly more aggressive as $\beta$ gets bigger, as will be studied in depth later.

The formulation readily extends to multi-stage problems, where decisions are made stage by stage as the scenario gradually reveals itself.
% such as structure, complexity and approximability, as well as on effective decomposition algorithms. 
Let $t=1, \cdots, T$ labels the stages sequentially, and let $x_t$ and $\zeta_t$ denote the decision and scenario for stage $t$.
% The decision variable $x$ now consists of $T$ subverctors $(x_1, \cdots, x_T)$, with the stage decision $x_t$ for stage $t$. 
% Likewise, a scenario breaks into stage scenarios: $\zeta = (\zeta_1, \cdots, \zeta_T)$. 
To standardize, a stage decision $x_t$ is always carried out before the stage scenario $\zeta_t$ is realized and known to the decision maker, which is without loss of generality: 
If a stage scenario is realized before any decision is made, a dummy decision can be added in the very beginning to make a standardized formulation.
Let $x=(x_1,\cdots,x_T)$ and $\zeta=(\zeta_1,\cdots,\zeta_T)$ for the entire decision and scenario.

It is assumed that decision makers neither know nor influence future stage scenarios by their actions, which is known as {\em nonanticipativity} in stochastic programming. % that the realization of scenarios is independent of decisions.
	%multi-stage stochastic programming (MSP) is a framework for 
	%sequential decision making under uncertainty where the decision 
	%space is typically high dimensional and involves complicated 
	%constraints, and the uncertainty is modeled by a general stochastic 
	%process. In the traditional risk neutral setting, 
	%the goal is to find a sequence of decisions or a policy so as to 
	%optimize an expected value objective. 
%The multi-stage scenarios in $\mathcal{U}$ can be organized into a scenario tree, with stage scenarios $\zeta_t$ as nodes so that a path from the root (the state before the first stage) to a leaf node corresponds to a scenario $(\zeta_1, \cdots, \zeta_T)$.
%As a stage scenario realizes itself, one moves along a particular path in the scenario tree. 
%The sub-tree from a node contains all possible future scenarios from that node.
Specifically, when making decision $x_t$ for stage $t$, only the partial scenario $\zeta_{1:t}= (\zeta_1, \cdots, \zeta_{t-1})$ is known, reducing the uncertainty set to $\mathcal{U}(\zeta_{1:t}) = \{ \omega \in \mathcal{U} : \omega_{1:t} = \zeta_{1:t} \}$.
With reduced $\mathcal{U}(\zeta_{1:t})$, the feasible set $X_{\mathcal{U}(\zeta_{1:t})}=\bigcap_{\omega \in \mathcal{U}(\zeta_{1:t})} X_\omega$ expands: $X_{\mathcal{U}(\zeta_{1:t})}\subseteq X_{\mathcal{U}(\zeta_{1:t+1})}$ as $\mathcal{U}(\zeta_{1:t+1})\subseteq\mathcal{U}(\zeta_{1:t})$.
Let $x_{1:t} = (x_1, \cdots, x_{t-1})$ be all stage decisions before stage $t$, and $h_t=(x_{1:t}, \zeta_{1:t})$ be the history. 
%By nonanticipativity, the stage decision $x_t$ depends only on $h_t$. 
The feasible set compatible with $h_t$ is $X(h_t)=\{y \in X_{\mathcal{U}(\zeta_{1:t})}: y_{1:t}=x_{1:t}\}$.
Let $X_t(h_t)=\{y_t : y \in X(h_t)\}$ be the stage feasible set given $h_t$, and $\mathcal{U}_t(\zeta_{1:t}) = \{ \omega_t : \omega \in \mathcal{U}(\zeta_{1:t}) \}$ be the  stage scenario set given $\zeta_{1:t}$.
A history $h_t$ is feasible if and only if each stage decision $x_\tau$ in $h_t$ is feasible given the sub-history $h_\tau$ in $h_t$ before stage $\tau$:
\begin{equation}\label{defn:feasible-history}
x_\tau\in X_\tau(h_\tau), \tau = 1, \cdots, t-1. % \text{ where } h_\tau=(x_{1:\tau}, \zeta_{1:\tau}). 
\end{equation}
Let $H_t$ denote the set of all feasible histories before stage $t$.

	% may use parameter p_t = p_t(x_{1:t-1}, \zeta_{1:t-1})
	% to formulate the case where action can influence parameters.
	%TODO: combine with a standard MP to illustrate $\mathcal{U}$, and $X$.
	
The ARM criterion for multi-stage problems can be defined recursively.
Let $r(x, \zeta)$ denote the total reward over all stages, which can be either accrued over stages or received all at once in the end.
Let $r^*(\zeta) = \max_{x \in X_\zeta} r(x,\zeta)$ be the ex post optimal. 
%where $X(\zeta) = \{ (x_1, \cdots, x_T): x_t \in X_t(h_t), t=1, \cdots, T\}$ is the set of all actions compatible with $\zeta$.  
% Note: 1. X(\zeta) is not needed to define r^*(\zeta), as \zeta is known before-hand, x can be chosen from X_\zeta. 2. there is $X(\zeta) \in X_\zeta$. 
The regret for a completed history $h_{T+1} = (x, \zeta)$ is  
\begin{equation}\label{defn:D[T](h)}
  D_T(h_{T+1}; \beta)= \beta r^*(\zeta) - r(x, \zeta).
\end{equation}
Work from $D_T(h_{T+1}; \beta)$, the regret guarantees of $\bar{D}_t(x_t, h_t; \beta)$ for $x_t, h_t$ and $D_{t-1}(h_t;\beta)$ for $h_t$ are found by {\em backward induction} for $t=T, \cdots, 1$: 
\begin{eqnarray}
	\bar{D}_t(x_t, h_t; \beta) &=& 
	\max_{\zeta_t \in \mathcal{U}_t(\zeta_{1:t})} D_t(h_{t+1}; \beta), 
	\label{defn:D[t](h)} \\
	D_{t-1}(h_t;\beta) &=& \min_{x_t\in X_t(h_t)} \bar{D}_t(x_t, h_t;\beta),\nonumber\\
	&=& \min_{x_t\in X_t(h_t)} \max_{\zeta_t \in \mathcal{U}_t(\zeta_{1:t})} D_t(h_{t+1};\beta). \label{defn:Dt(h)} 
  \end{eqnarray}
% \begin{equation}\label{defn:barD[t](h)}
%   \bar{D}_t(x_t, h_t; \beta) = 
%   \max_{\zeta_t \in \mathcal{U}_t(\zeta_{1:t})} D_t(h_{t+1}; \beta),
% \end{equation}
where $h_{t+1}$ is formed by adding $x_t$ and $\zeta_t$ into $h_t$.
Note that \eqref{defn:D[t](h)} is often referred to as the ``adversarial problem'', as if an almighty adversary is always making the decision maker regret the most.
% An optimal stage action $x_t$ is chosen to minimize the regret guarantee
%The definition (\ref{defn:D[t](h)}) can be applied recursively for $t=T, \cdots, 1$ backwards, which 
As $h_1$ is empty, let $D(\beta) \equiv D_0(h_1;\beta)$ for the best regret guarantee, and the performance guarantee is still $\check{r}(r, \beta) = \beta r - D(\beta)$. % qualifying reward
This completes the {\em plain formulation}.

An alternative formulation is based on policies. 
A feasible policy $\pi$ is a sequence of functions $\pi = \{ \pi_t: \pi_t(h_t) \in X_t(h_t), \forall h_t \in H_t, t=1, \cdots, T \}$, which makes stage decisions by $x_t = \pi_t(h_t)$, with nonanticipativity already baked-in.
% Random policies can be defined by replacing $X_t(h_t)$ by a set of probabilistic distributions on $X_t(h_t)$,  but the focus is on deterministic policies. (This randomization gives behavior policies)
% Note: one may also randomize among deterministic policies (mixed policies)
% see https://en.wikipedia.org/wiki/Strategy_(game_theory)
% Behavior strategy
% While a mixed strategy assigns a probability distribution over pure strategies, a behavior strategy assigns at each information set a probability distribution over the set of possible actions. 
% Defined for any feasible history, it works well with backward induction.
Deterministic policies are the focus, though random ones are possible. % may be defined by replacing $X_t(h_t)$ with a set of probabilistic distributions on $X_t(h_t)$.
Apply a policy $\pi$ to a $\zeta\in U$ from stage 1 to T, the whole decision $x$ is uniquely determined, which can be simply denoted as $x=\pi(\zeta)$.
All partial histories before stage $t$ that is realizable under $\pi$ is defined as  $H_t^\pi=\{(x_{1:t},\zeta_{1:t}) : x=\pi(\zeta), \zeta \in \mathcal{U}\}$.
The restriction of $\pi_t$ to $H_t^\pi$ defines a $\ddot{\pi}_t$, which gives a pruned policy $\ddot{\pi}=\{\ddot{\pi}_t: t=1,\cdots,T\}$. % which is practically simpler and theorectically useful.
Let $\Pi$ be the set of all feasible policies, and let $\ddot{\Pi}=\{\ddot{\pi}: \pi \in \Pi\}$ be the pruned set. %, where the double-dot accent indicates the pruning operation.
%Note that if a policy $\pi$ is strictly followed, 
%then $h^x_t$ can be recursively computed from $h^\zeta_t$ by
%$x_i=\pi_i(h_i)=\pi_i(h_i)$ for $i=1, 2, \cdots, t-1$,
%or shorthanded as $x_{1:t-1} = \pi(\zeta_{1:t-1})$.
%In this case, the policy is given without $x_{1:t-1}$,
%$\pi_t(\zeta_{1:t-1}) = \pi_t(\pi(\zeta_{1:t-1}), \zeta_{1:t-1})$.
    %It is possible to start applying a policy from stage $t$ on with an arbitrary history $h_t$.
	% but if a policy $\pi$ is applied throughout the stages, then $x_{1:t}$ and $x_t$ can be computed from $\zeta_{1:t}$, it is succinct to have $x^\pi_{1:t+1}(\zeta_{1:t}) = (x^\pi_1, \cdots, x^\pi_t)$.

The regret guarantee of $\pi$ is also found by backward induction.
For a fully developed history $h_{T+1}=(x, \zeta)$, the regret is still
\begin{equation}\label{defn:D[T](pi)}
  D_T^\pi(h_{T+1};\beta) = \beta r^*(\zeta) - r(x, \zeta).
\end{equation}
The regret guarantee given $h_t=(x_{1:t},\zeta_{1:t})$ is defined recursively by
\begin{equation}\label{defn:D[t](pi)}
  D_{t-1}^\pi(h_t;\beta) = 
  \max_{\zeta_t \in \mathcal{U}_t(\zeta_{1:t})} 
  D_{t}^\pi (h^\pi_{t+1}; \beta), \; t = T, \cdots, 1,
\end{equation}
where $h^\pi_{t+1} = ((x_{1:t},\pi_t(h_t)), (\zeta_{1:t}, \zeta_t))$ evolves from $h_t$ under $\pi$. % if $\zeta_t$ is realized in stage $t$
Apply (\ref{defn:D[t](pi)}) recursively for the policy's regret guarantee $D^\pi(\beta) \equiv D_{0}^\pi(h_1, \beta)$:
\begin{eqnarray}\label{defn:D[0](pi)}
  D^\pi(\beta) &=& \max_{\zeta_{1} \in \mathcal{U}_1(\zeta_{1:1})} D_{1}^\pi(h^\pi_2; \beta) \nonumber\\
  &=&\max_{\zeta_{1} \in \mathcal{U}_1(\zeta_{1:1})} \max_{\zeta_{2} \in \mathcal{U}_2(\zeta_{1:2})} D_{2}^\pi(h^\pi_3; \beta)\nonumber\\
  &=&\max_{\zeta_{1} \in \mathcal{U}_1(\zeta_{1:1})} \cdots \max_{\zeta_{T} \in \mathcal{U}_T(\zeta_{1:T})} D_{T}^\pi (h^\pi_{T+1}; \beta)\nonumber\\
  &=&\max_{\zeta \in \mathcal{U}} D_{T}^\pi(h^\pi_{T+1}; \beta)
\end{eqnarray}
where $h^\pi_{T+1}=(\pi(\zeta),\zeta)$.
The {\em policy-based formulation} chooses a policy to minimize the regret guarantee:
\begin{eqnarray}\label{defn:policy-form}
  \min_{\pi\in\Pi} D^\pi(\beta) = 
  \min_{\pi\in\Pi} \max_{\zeta\in\mathcal{U}} \beta r^*(\zeta) - r(\pi,\zeta),
\end{eqnarray}
where $r(\pi,\zeta)=r(\pi(\zeta), \zeta)$, as $\pi$ in \eqref{defn:policy-form} plays the role of $x$ in \eqref{defn:arm-beta}.

Some comments on the policy-based formulation. %conceals the multi-stage process but
First, all stage decisions are made by  $\pi$ in \eqref{defn:policy-form}, producing only realizable histories $H_t^\pi$, thus  $\Pi$ can be replaced by $\ddot{\Pi}$. % in \eqref{defn:policy-form}.
% Hence for an unrealizable history $h_t$ under an optimal policy $\pi^*$, the guarantee $D_{t-1}^{\pi^*}(h_t;\beta)$ is unsolvable by \eqref{defn:policy-form}, which invokes
%Clearly, it is fine to replace $\Pi$ by $\ddot{\Pi}$ in \eqref{defn:policy-form}, and the policy $\ddot{\pi}^*$ pruned from $\pi^*$ keeps the minimum guarantee for all realizable histories.
Second, for problems with states, any $h_t\in H_t$ maps to a state $s(h_t)$, which may replace $h_t$ to reduce complexity with fewer states than histories.
Finally, there could be many optimal policies from \eqref{defn:policy-form} that may differ in value for some $D_{t-1}^{\pi}(h_t;\beta)$. 
The set of optimal policies can be reduced by the principle of optimality (\citealt{bellman1954theory}), which requires optimality recursively: 
%that the subpolicies of an optimal policy are themselves optimal: 
A policy $\pi^*$ is optimal if and only if $D^{\pi^*}_{t-1}(h_t; \beta ) \le D^{\pi}_{t-1}(h_t; \beta)$ for all $\pi\in \Pi$ and $h_t \in H_t, t=1,\cdots,T$. 

\section{Theoretical Analysis.} \label{sec:prop}
% Fundamental properties of the ARM criterion are studied in this section.
In this section, formulation equivalence and convexity preservation are established first, then a new approach to competitive ratio analysis is proposed on a solid theoretical basis, and finally a framework  for fine control of conservatism is studied.

	%which lays the theoretic foundations for
	%its applications in RO. 
	%For the popular competitive ratio analysis, 
	%the ARM enables a novel approach that can greatly reduce 
	%its complexity of analysis to that of the minimax regret analysis.

%The correspondence and equivalence between the formulations is first investigated.

%Backward induction and sub-game perfect equilibrium (zero-sum game) policy,
%the policy-based formulation.
 
\begin{lem}%[Formulation Equivalence] 
\label{lem:optimality}
Under the principal of optimality, the plain and the policy-based formulations are equivalent in that any optimal policy $\pi^*\in\Pi$ satisfies
%there is $\Pi^* = \Pi^*_1 = \Pi^*_2$, where $\Pi^*\subseteq\Pi$ is the set of optimal policies for \eqref{defn:policy-form}, and
% $\Pi^*_1 = \{ \pi \in \Pi: D_{t-1}(h_{t}; \beta) = D_{t-1}^{\pi}( h_{t}; \beta), \forall h_t \in H_t, t=1, \cdots, T+1 \}, \Pi^*_2 =\{\pi \in \Pi: \pi_t(h_t) \in \argmin_{x_t\in X_t(h_t)}\; \max_{\zeta_t \in \mathcal{U}_t(\zeta_{1:t})}\; D_t(h_{t+1};\beta), \forall h_t \in H_t, t=1,\cdots, T\}$
\begin{eqnarray} % \label{eqn:value-equivalence}
	D_{t-1}(h_{t}; \beta) =	D_{t-1}^{\pi^*}( h_{t}; \beta), \forall h_t \in H_t, t=1, \cdots, T+1, \label{eqn:value-equivalence} \\ 
	\pi^*_t(h_t) \in \argmin_{x_t\in X_t(h_t)} \max_{\zeta_t \in \mathcal{U}_t(\zeta_{1:t})} D_t(h_{t+1};\beta), \forall h_t \in H_t, t=1, \cdots, T, \label{defn:optimal-policy}
\end{eqnarray}
where $\zeta_{1:t}$ belongs to $h_t=(x_{1:t},\zeta_{1:t})$.
% \begin{equation} \label{eqn:value-equivalence}
% 	D_{t-1}(h_{t}; \beta) = 
% 	D_{t-1}^{\pi^*}( h_{t}; \beta), 
% 	\forall h_t \in H_t, t=1, \cdots, T+1,
% \end{equation}
% (ii) if a policy $\pi^* \in \Pi$ satisfies \eqref{eqn:value-equivalence}, then it is optimal, 
% and (iii) if a policy $\pi^* \in \Pi$ is optimal, then it satisfies \eqref{defn:optimal-policy}, which loops back to (i).
\end{lem}

Please see \ref{app:optimality} for the proof of Lemma \ref{lem:optimality}. 
It is handy to have both formulations, as the plain one solves a problem stage by stage, while the policy-based one facilitates theoretical analysis. 
The interchangeability principle of \cite{shapiro2017interchangeability} can obtain a weaker result, as it does not observe the principle of optimality. 
%Conversely, the proof here can be adapted to make an alternative proof of the interchangeability principle.
Also note that there is $D^{\pi^*}(\beta)=D(\beta)$ by \eqref{eqn:value-equivalence}.

Convexity preservation is crucial for computational tractability, and the analysis employs the concepts of consistent convexity and convex combination of pruned policies.
Consistent convexity simply means that the problem $\max \{ r(x, \zeta): x \in X_\zeta \}$ is consistently convex for all $\zeta \in \mathcal{U}$. %with a convex $X_\zeta$ and an objective $r(x, \zeta)$ jointly concave in $x$. 
A convex combination of two pruned policies $\ddot{\pi}_1, \ddot{\pi}_2\in \ddot{\Pi}$ by a $\lambda \in [0,1]$ is defined as $\ddot{\pi}(\zeta)=\lambda \ddot{\pi}_1(\zeta)+(1-\lambda)\ddot{\pi}_2(\zeta)$ for all $\zeta\in \mathcal{U}$, which can be written as $\ddot{\pi} = \lambda \ddot{\pi}_1 + (1-\lambda)\ddot{\pi}_2$ with $\ddot{\pi}$ regarded as a vector made of $\ddot{\pi}(\zeta), \forall\zeta\in \mathcal{U}$.

% \begin{lem} %[Continuity] 
% \label{lem:pruned-convexity}
% With consistent convexity, the pruned feasible set $\ddot{\Pi}$ is convex.
% \end{lem}
% \begin{proof} 
% 	For $\ddot{\pi}_1, \ddot{\pi}_2\in \ddot{\Pi}$, the history $h^{\ddot{\pi}_i}_{T+1}=(\ddot{\pi}(\zeta),\zeta)$ must be feasible and satisfy \eqref{defn:feasible-history} for all $\zeta\in \mathcal{U}$.
% 	For $\ddot{\Pi}$ to be convex, it suffices to show that policy $\ddot{\pi} = \lambda \ddot{\pi}_1 + (1-\lambda)\ddot{\pi}_2$ combined by $\lambda \in [0,1]$ is feasible.
% 	Apply both $\ddot{\pi}_i$ to an arbitrary $\zeta$ to have $x^i=\ddot{\pi}_i(\zeta)$, clearly $h^{\ddot{\pi}_i}_{T+1}=(x^i,\zeta)$ satisfy \eqref{defn:feasible-history}, so that there exists $y^{it} \in X_{\mathcal{U}(\zeta_{1:t})}$ such that $y^{it}_{1:t} = x^i_{1:t}$ for $i=1,2$ and $t=1,\cdots, T$.
% 	Consistent convexity implies that any $X_{\mathcal{U}(\zeta_{1:t})}$ is convex, which ensures that $y^{t} = \lambda y^{1t} + (1-\lambda) y^{2t} \in X_{\mathcal{U}(\zeta_{1:t})}$, hence for $x=\ddot{\pi}(\zeta), h^{\ddot{\pi}}_{T+1}=(x,\zeta)$ there is $x_t \in X_t(h^{\ddot{\pi}}_t)$ for $t=1,\cdots,T$ due to $y^{t}$, so $h^{\ddot{\pi}}_{T+1}$ satisfies \eqref{defn:feasible-history} and $\ddot{\pi}\in \ddot{\Pi}$.
% \end{proof}

\begin{thm}%[Convexity Preservation]
	\label{thm:convexity-preservation}
	With consistent convexity, the policy-based formulation \eqref{defn:policy-form} restricted to pruned policies $\ddot{\Pi}$ is a convex problem, and all subproblems \eqref{defn:Dt(h)} in the plain formulation are convex problems. 
	%with convex sets $\ddot{\Pi}$ and $X_t(h_t)$, and convex objectives $D^{\ddot{\pi}}(\beta)$ and $\bar{D}_t(x_t, h_t; \beta)$.
	% \begin{itemize}
	% 	\item[1.] Both the pruned feasible set $\ddot{\Pi}$ and $X_t(h_t)$ are convex sets;
	% 	\item[2.] Both the objectives $D^{\ddot{\pi}}(\beta)$ and $\bar{D}_t(x_t, h_t; \beta)$ are convex functions.
	% \end{itemize}
\end{thm}
\begin{proof}
	First prove the convexity of the domains.
	As a projection of the convex set $X(h_t)$, the convexity of $X_t(h_t)$ is obvious,	but the convexity of $\ddot{\Pi}$ needs some explanation. 
	For $\ddot{\pi}^1, \ddot{\pi}^2\in \ddot{\Pi}$, the histories $h^{\ddot{\pi}^i}_{T+1}=(\ddot{\pi}(\zeta),\zeta), i=1,2$ satisfy the feasibility condition \eqref{defn:feasible-history} for all $\zeta\in \mathcal{U}$.
	For $\ddot{\Pi}$ to be convex, it must be shown that policy $\ddot{\pi} = \lambda \ddot{\pi}^1 + (1-\lambda)\ddot{\pi}^2$ for any $\lambda \in [0,1]$ also satisfies \eqref{defn:feasible-history} to have $\ddot{\pi}\in \ddot{\Pi}$.
	Apply both $\ddot{\pi}^i, i=1,2$ to an arbitrary $\zeta$ to have $x^i=\ddot{\pi}^i(\zeta), i=1,2$. As $h^{\ddot{\pi}^i}_{T+1}=(x^i,\zeta), i=1,2$ satisfy \eqref{defn:feasible-history}, there exists $y^{it} \in X_{\mathcal{U}(\zeta_{1:t})}$ such that $y^{it}_{1:t} = x^i_{1:t}$ for $i=1,2$ and $t=1,\cdots, T$.
	Consistent convexity implies that any $X_{\mathcal{U}(\zeta_{1:t})}$ is convex, which ensures $y^{t} = \lambda y^{1t} + (1-\lambda) y^{2t} \in X_{\mathcal{U}(\zeta_{1:t})}$, hence for $x=\ddot{\pi}(\zeta), h^{\ddot{\pi}}_{T+1}=(x,\zeta)$ there is $x_t \in X_t(h^{\ddot{\pi}}_t)$ for $t=1,\cdots,T$, so $h^{\ddot{\pi}}_{T+1}$ satisfies \eqref{defn:feasible-history} and there is $\ddot{\pi}\in \ddot{\Pi}$.

	% For $\ddot{\pi}^1, \ddot{\pi}^2\in \ddot{\Pi}$, the histories $h^{\ddot{\pi}^i}_{T+1}=(\ddot{\pi}(\zeta),\zeta), i=1,2$ satisfy the feasibility condition \eqref{defn:feasible-history} for all $\zeta\in \mathcal{U}$, thus for .
	% Let policy $\ddot{\pi} = \lambda \ddot{\pi}^1 + (1-\lambda)\ddot{\pi}^2$ with $\lambda \in [0,1]$.
	% Consistent convexity implies that  $X_{\mathcal{U}(\zeta_{1:t})}$ is convex, and so is $X(h_t)$ as a slice of $X_{\mathcal{U}(\zeta_{1:t})}$, which makes $X_t(h_t)$ convex too as it is a projection of $X(h_t)$. 
	% Therefore, policy $\ddot{\pi} = \lambda \ddot{\pi}_1 + (1-\lambda)\ddot{\pi}_2$ with $\lambda \in [0,1]$ also satisfies \eqref{defn:feasible-history}, so there is $\ddot{\pi}\in \ddot{\Pi}$ and $\ddot{\Pi}$ is convex.

	Next show the convexity of the objectives.
	The regret of any $\pi$ with a given $\zeta\in\mathcal{U}$ is $E_\zeta(\pi)\equiv\beta r^*(\zeta) - r(\pi, \zeta)$. %E is a function of \pi
	For $\ddot{\pi}_1, \ddot{\pi}_2\in \ddot{\Pi}$ and a $\lambda\in[0,1]$, let $\ddot{\pi} = \lambda \ddot{\pi}_1 + (1-\lambda)\ddot{\pi}_2 \in \ddot{\Pi}$. 
	For any $\zeta\in\mathcal{U}$, there is $\ddot{\pi}(\zeta)=\lambda \ddot{\pi}_1(\zeta) + (1-\lambda)\ddot{\pi}_2(\zeta)$, thus $E_\zeta(\ddot{\pi}) \leq \lambda E_\zeta(\ddot{\pi}_1) + (1-\lambda) E_\zeta(\ddot{\pi}_2)$ as $r(x,\zeta)$ is concave by consistent convexity, so $E_\zeta(\pi)$ is convex on $\ddot{\Pi}$.
	The objective at $\ddot{\pi}\in \ddot{\Pi}$ in \eqref{defn:policy-form} is 
	$D^{\ddot{\pi}}(\beta)=\max_{\zeta\in\mathcal{U}} E_\zeta(\ddot{\pi})$,	which is convex in $\ddot{\pi}$ as a pointwise max of convex functions defined on the domain of $\ddot{\Pi}$.
	
	It remains to show that $\bar{D}_t(x_t, h_t; \beta)$ in \eqref{defn:Dt(h)} is convex in $x_t$ for any $h_t=(x_{1:t},\zeta_{1:t})\in H_t$. 
	%X_{1:t+1}(\zeta_{1:t}) \equiv \{ y_{1:t+1}: (y_{1:t}, \zeta_{1:t})\in H_{t}, y_t \in X_t(y_{1:t}, \zeta_{1:t})\} \equiv \{ \ddot{\pi}_{1:t+1}(\zeta_{1:t}): \ddot{\pi} \in \ddot{\Pi} \}$. % this
	Let $\ddot{\Pi}(h_t) = \{ \ddot{\pi}\in \ddot{\Pi}:\ddot{\pi}_\tau(\zeta_{1:\tau})=x_\tau, \tau = 1, \cdots, t-1 \}$ and $\ddot{\Pi}(x_t, h_t) =\{ \ddot{\pi}\in \ddot{\Pi}(h_t): \ddot{\pi}_t(\zeta_{1:t})=x_t\}$, so that $\ddot{\Pi}(h_t) =\bigcup_{x_t\in X_t(h_t)} \ddot{\Pi}(x_t, h_t)$.
	Due to nonanticipativity, any $\ddot{\pi}\in \ddot{\Pi}(h_t)$ satisfies $\ddot{\pi}_t(\zeta'_1) = \ddot{\pi}_t(\zeta'_2) $ for any $\zeta'_1, \zeta'_2 \in \mathcal{U}(\zeta_{1:t})$, thus $\ddot{\pi}$ can be represented as a vector with only one $x_t$ component instead of many copies of the same $x_t$ for each $\zeta'\in \mathcal{U}(\zeta_{1:t})$.
	Both $\ddot{\Pi}(x_t, h_t)$ and $\ddot{\Pi}(h_t)$ are convex as slices of $\ddot{\Pi}$, and $\ddot{\Pi}(x_t, h_t)$ is a slice of $\ddot{\Pi}(h_t)$.
	Apply Lemma \ref{lem:optimality} to the subproblem \eqref{defn:D[t](h)} as an independent problem with an uncertainty set $\mathcal{U}(\zeta_{1:t})$ and a dummy decision $x_t$ so that $\ddot{\Pi}(x_t, h_t)$ has all pruned policies for it, to have $\bar{D}_t(x_t, h_t; \beta) = \min_{\ddot{\pi}\in \ddot{\Pi}(x_t, h_t)} g(\ddot{\pi})$, where $g(\ddot{\pi}) \equiv \max_{\zeta \in \mathcal{U}(\zeta_{1:t})} E_\zeta(\ddot{\pi})$ is convex over $\ddot{\Pi}(x_t, h_t)$.
	The epigraph $\mathbf{epi }\; \bar{D}_t(X_t(h_t), h_t; \beta) = \{(x_t, v): (\ddot{\pi}, v) \in \mathbf{epi }\; g(\ddot{\Pi}(h_t)) \text{ for some } \ddot{\pi} \in \ddot{\Pi}(x_t, h_t) \}$ is convex as a projection of the convex set $\mathbf{epi }\; g(\ddot{\Pi}(h_t)) =\{ (\ddot{\pi}, v): v \ge g(\ddot{\pi}) \} $ onto $x_t$ as a component of $\ddot{\pi}$. % see page 88 of "convex optimization" by Boyd.
	Thus the objective $\bar{D}_t(x_t, h_t; \beta)$ in \eqref{defn:Dt(h)} is convex in $x_t$.
\end{proof}	% TODO: how to exploit this convexity?

Convexity preservation not only facilitates theoretical analysis, but ensures global convergence of numerical algorithms such as robust policy/value iteration (see \citealt{moos2022robust} for a survey).
Note that the proof does not get into the ``adversarial problem'' of \eqref{defn:D[t](h)}, which may be challenging to deal with when looking for analytical solutions or designing numerical algorithms.
By the way, Theorem \ref{thm:convexity-preservation} also applies to the three commonly used criteria as special cases of the ARM criterion, whose tractability also comes from convexity preservation.
The ARM criterion resembles the minimax regret criterion, where the only difference is in the value of $\beta$, thus a reasonable conjecture is that they give formulations of similar tractability.

\subsection{Competitive Ratio.}
Competitive ratio has been applied when relative regret is more appropriate than absolute regret, and analytical solutions are derived sometimes even with discrete variables (e.g. \citealt{wang2022robust}).
Among equivalent variants, the competitive ratio for reward maximization can be defined as 
	\begin{eqnarray}\label{defn:cratio}
		\gamma^* = \max_{\pi\in\Pi} \min_{\zeta\in\mathcal{U}} {r_\zeta(\pi) / r^*(\zeta)},
	\end{eqnarray}
where $r^*(\zeta) > 0$ for all $\zeta\in\mathcal{U}$ is assumed in general. 
%To pave a new way to competitive ratio analysis, the basic property of continuity is studied first.

\begin{lem}%[Continuity] 
\label{lem:continuity}
  The regret guarantee $D_{t-1}(h_t; \beta)$ for $t=1, \cdots, T+1$ with an arbitrary history $h_t\in H_t$ is continuous in $\beta$.
\end{lem}
\begin{proof} Use backward induction on $t$.
When $t=T+1$, it is clear that $D_T(h_{T+1}; \beta)$ is continuous in $\beta$ according to \eqref{defn:D[T](h)}, which completes the initial step. 
The induction step assumes $D_{t}(h_{t+1}; \beta)$ is continuous in $\beta$, then shows the same for $D_{t-1}(h_t; \beta)$. 
It is clear that $\bar{D}_t(x_t, h_t; \beta)$ is continuous in $\beta$ as it is a point-wise max of continuous functions by \eqref{defn:D[t](h)}. 
Likewise, $D_{t-1}(h_t; \beta)$ is also continuous with regard to $\beta$ by \eqref{defn:D[t](h)}.
\end{proof}

%The continuity of Lemma \ref{lem:continuity} is a basic property useful for other analysis later, such as in Lemma \ref{lem:unique-cratio} for the existence and uniqueness of CR. 

The continuity of $D(\beta)$ suggests that the optimal policy and the conservatism may change in a ``smooth'' manner with $\beta$, which is observed in the one-way trading problem studied later. % added by the suggestion of reviewer #1 major comment #16.
Next the slope bounds for $D(\beta)$ is given, which establishes that $D(\beta)$ strictly increases if $r^*(\zeta) > 0$ for all $\zeta\in\mathcal{U}$. 

\begin{thm}%[Slope Bounds] 
\label{thm:slope-bounds}
	For $\beta_1 < \beta_2$, let $\pi^*_i, i\in\{1,2\}$ be an optimal policy for $\beta = \beta_i$, and $\zeta^*_{ij} = \argmax_{\zeta\in\mathcal{U}} \beta_i r^*(\zeta) - r(\pi^*_j,\zeta), i,j\in\{1,2\}$, then there is 
		\begin{eqnarray}\label{ineq:monotony}
		r^*(\zeta^*_{21}) \ge {D(\beta_2) - D(\beta_1)\over \beta_2-\beta_1 } \ge r^*(\zeta^*_{12}).
		\end{eqnarray}
\end{thm}
\begin{proof} 
By the definition of $\pi^*_2$ and $\zeta^*_{12}$, as well as Lemma \ref{lem:optimality}, there is
\begin{eqnarray*}
		D(\beta_1) &=& \min_{\pi\in\Pi} \max_{\zeta\in\mathcal{U}} 
		\beta_1 r^*(\zeta) - r(\pi,\zeta)\\
		&\le& \max_{\zeta\in\mathcal{U}} \beta_1 r^*(\zeta) - r(\pi^*_2,\zeta) \\
		&=& \beta_1 r^*(\zeta^*_{12}) - r(\pi^*_2,\zeta^*_{12}).
\end{eqnarray*} 
And there is $D(\beta_2) = \max_{\zeta\in\mathcal{U}} 
		\beta_2 r^*(\zeta) - r(\pi^*_2,\zeta) \ge \beta_2 r^*(\zeta^*_{12}) - r(\pi^*_2,\zeta^*_{12})$. 
Therefore $D(\beta_2) - D(\beta_1) \ge
		(\beta_2 - \beta_1)\ r^*(\zeta^*_{12})$. 
Similarly,
		\begin{eqnarray*}
			D(\beta_2) &=& \min_{\pi\in\Pi} \max_{\zeta\in\mathcal{U}} 
			\beta_2 r^*(\zeta) - r(\pi,\zeta)\\
			&\le& \max_{\zeta\in\mathcal{U}} \beta_2 r^*(\zeta) - r(\pi^*_1,\zeta) \\
			&=& \beta_2 r^*(\zeta^*_{21}) - r(\pi^*_1,\zeta^*_{21}).
		\end{eqnarray*} 
And there is 
		$D(\beta_1) = \max_{\zeta\in\mathcal{U}} 
		\beta_1 r^*(\zeta) - r(\pi^*_1,\zeta) \ge \beta_1 r^*(\zeta^*_{21}) - r(\pi^*_1,\zeta^*_{21})$.
Thus $D(\beta_2) - D(\beta_1) \le
		(\beta_2 - \beta_1)\ r^*(\zeta^*_{21})$.  
Therefore, (\ref{ineq:monotony}) follows immediately.
\end{proof}

% For example, clearing all holdings is a way to avoid loss in the stock market. 
% In applications when this assumption is not satisfied,  a reference performance can be provided by a function $r'(\zeta)$ that satisfies $r'(\zeta) \le r^*(\zeta), \zeta\in\mathcal{U}$, so that the reduced loss above the reference $r_+(x,\zeta) = r(x,\zeta) - r'(\zeta)\ge 0$ may be employed instead as the new reward function. 
% Such function $r'(\zeta)$ clearly exists, for example,  $f(\zeta) = \min_{x\in X} r(x,\zeta)$ or simply $f(\zeta) = \min_{\zeta\in\mathcal{U}} r^*(\zeta)$, to make $r_+(x,\zeta)$ satisfying the loss prevention assumption. 
% which is assumed by default in this paper.

With continuity and monotonicity from Lemma \ref{lem:continuity} and Theorem \ref{thm:slope-bounds}, it is ready to lay the theoretical foundation for a new approach to competitive ratio analysis.

% Note that r*(\zeta_i)>0 for all \zeta \in U does not imply D(0)<0: 
% a counter example can be constructed: 
% U={\zeta_i = 2i\pi/3, i=0,1,2}, x\in[0,2\pi] with r(x,\zeta)=sin(x+\zeta).
% there is \max_{x} \min_{\zeta\in U} r(x,\zeta) < 0. 
% which means D(0)>0 while r*(\zeta)>0 for all \zeta\in U.

\begin{thm}%[ARM for CR] 
\label{thm:cr-equivalence}
With $r^*(\zeta) > 0$ for all $\zeta\in\mathcal{U}$, the competitive ratio $\gamma^*$ is the unique solution to $D(\beta)=0$, and problem \eqref{defn:policy-form} with $\beta=\gamma^*$ has the same optimal policies as problem \eqref{defn:cratio}.
%moreover, the sufficient and necessary condition for a unique and nonnegative competitive ratio is $D(0)\le 0$.
% and $r(x,\zeta)$ is bounded below, then $D(\beta)=0$ has a unique solution $\beta_0\le 1$ that equals the competitive ratio $\gamma^*$, 
\end{thm}
\begin{proof} 
Positive $r^*(\zeta)$ ensures a strictly increasing $D(\beta)$ by Theorem \ref{thm:slope-bounds}, thus $D(\beta)=0$ can never have more than one solution. 
%with $r(x,\zeta)$ bounded below there is $\lim_{\beta\downarrow-\infty} D(\beta) = -\infty$, and $D(1)\ge 0$ is obvious, $\beta_0 \le 1$ by continuity from Lemma \ref{lem:continuity}.
% It needs to show that any optimal $\pi^*$ for (\ref{defn:policy-form}) also solves (\ref{defn:cratio}), and vice versa.
Start from (\ref{defn:cratio}) and proceed as follows:
{\allowdisplaybreaks\begin{eqnarray*}
	&&\left\{\begin{array}{l}\displaystyle
		\gamma^* =  \max_{\pi\in\Pi} \min_{\zeta\in\mathcal{U}} r(\pi,\zeta)/r^*(\zeta)\\
		\displaystyle
		\pi^* \in \argmax_{\pi\in\Pi} \min_{\zeta\in\mathcal{U}} r(\pi,\zeta)/r^*(\zeta)
	\end{array}\right.\\
	&\Leftrightarrow&\left\{\begin{array}{l}\displaystyle
		\gamma^* =  \min_{\zeta\in\mathcal{U}} r(\pi^*,\zeta)/r^*(\zeta) \\
		\displaystyle
		\forall \pi\in\Pi: \gamma^* \geq \min_{\zeta\in\mathcal{U}} r(\pi,\zeta)/r^*(\zeta)
	\end{array}\right.\\
	&\Leftrightarrow&\left\{\begin{array}{l}
		\exists \zeta\in\mathcal{U}: \gamma^* = r(\pi^*,\zeta)/r^*(\zeta)\\
		\forall \zeta\in\mathcal{U}: \gamma^* \leq r(\pi^*,\zeta)/r^*(\zeta)\\
		\forall \pi \in \Pi, \exists \zeta\in\mathcal{U}: \gamma^* \geq r(\pi,\zeta)/r^*(\zeta)
	\end{array}\right.\\
	&\Leftrightarrow&\left\{\begin{array}{l}\displaystyle
		\exists \zeta\in\mathcal{U}: 0 = \gamma^* r^*(\zeta) - r(\pi^*,\zeta)\\
		\forall \zeta\in\mathcal{U}: 0 \geq \gamma^* r^*(\zeta) - r(\pi^*,\zeta)\\
		\forall \pi\in\Pi, \exists\zeta\in\mathcal{U}: 0 \leq \gamma^* r^*(\zeta) - r(\pi,\zeta)
	\end{array}\right.\\
	&\Leftrightarrow& \left\{\begin{array}{l}\displaystyle
		0=\max_{\zeta\in\mathcal{U}} \gamma^* r^*(\zeta) - r(\pi^*,\zeta)\\
		\displaystyle
		\forall \pi\in\Pi: 0\leq\max_{\zeta\in\mathcal{U}} \gamma^* r^*(\zeta) - r(\pi,\zeta)
	\end{array}\right.\\
	&\Leftrightarrow&\left\{\begin{array}{l}
		\displaystyle
		0 = \min_{\pi\in\Pi}\max_{\zeta\in\mathcal{U}} \gamma^* r^*(\zeta) - r(\pi,\zeta)\\
		\displaystyle
		\pi^*\in \argmin_{\pi\in\Pi}\max_{\zeta\in\mathcal{U}} \gamma^* r^*(\zeta) - r(\pi,\zeta)\\
	\end{array}\right.
\end{eqnarray*}}
As it can go both ways, the theorem is established.
% For the second part, since $D(\beta)$ is continuous and strictly increases due to $r^*(\zeta) > 0$ for all $\zeta\in\mathcal{U}$ and Lemma \ref{lem:continuity} and Theorem \ref{thm:slope-bounds}, and $D(1)\ge 0$ is obvious, thus $D(0)\le 0$ is sufficient for a unique $\gamma^*\in [0,1]$ that solves $D(\beta)=0$. 
% Conversely, if $D(0) > 0$, then $D(\beta) > 0$ for all $\beta>0$, which means $D(\beta)=0$ has no nonnegative solution.
\end{proof}  

So the ARM criterion recovers the relative regret criterion if $\beta$ is set to the competitive ratio. Note that \cite{averbakh2005computing} presents a similar result for single-stage problems, which is now extended to multi-stage problems. % so that more efficient numerical algorithms to compute competitive ratios are also possible based for multi-stage problems.
%, which could be negative for some problems. 
By Theorem \ref{thm:cr-equivalence}, a condition for a positive competitive ratio is derived:
\begin{lem} %[Uniqueness of CR]
\label{lem:unique-cratio}
	If $D(0)<0$ then $r^*(\zeta) > 0$ for all $\zeta\in\mathcal{U}$ and $\gamma^*\in(0,1]$.
\end{lem}
\begin{proof} 
Note that at $\beta = 0$ it gives the maximin reward criterion:
	\begin{equation*}
		D(0) = \min_{\pi\in\Pi}\max_{\zeta \in \mathcal{U}} - r(\pi,\zeta)
		= - \max_{\pi\in\Pi}\min_{\zeta \in \mathcal{U}}  r(\pi,\zeta).
	\end{equation*} 
 For an arbitrary $\bar{\zeta}$, there is
	\[
		%	\min_{\zeta\in\mathcal{U}} r^*(\zeta) = \min_{\zeta\in\mathcal{U}} \max_{\pi\in\Pi} r(x^\pi(\zeta),\zeta) 
		0 < - D(0)  = \max_{\pi\in\Pi} \min_{\zeta\in\mathcal{U}} r(\pi,\zeta) 
		\le \max_{\pi\in\Pi} r(\pi,\bar{\zeta}) \le r^*(\bar{\zeta}),
	\]
	Thus $D(\beta)$ strictly increases in $\beta$ by Theorem \ref{thm:slope-bounds}.
	Note that at $\beta=1$ it is the minimax regret criterion with $D(1) \ge 0 > D(0)$. So $D(\beta)=0$ has a unique solution $\beta_0$ that satisfies $\beta_0\in (0,1]$ by monotonicity and continuity of $D(\beta)$ from Lemma \ref{lem:continuity}.	
	The conclusion follows from Theorem \ref{thm:cr-equivalence}.
%ref: Boyd, Stephen; Vandenberghe, Lieven (2004), Convex Optimization, Cambridge University Press.
\end{proof}

A new approach to competitive ratio analysis comes straight out of Theorem \ref{thm:cr-equivalence}. 
If an analytical expression for $D(\beta)$ exists, then the competitive ratio can be found by simply solving $D(\beta)=0$.
This approach is generally simpler than directly dealing with the ratio in \eqref{defn:cratio}, as illustrated by the one-way trading problem later.
If $D(\beta)$ can be evaluated numerically, then the competitive ratio can be found by numerically solving $D(\beta)=0$ with modest extra computation cost, as in \cite{averbakh2005computing} for single-stage problems.
% \cite{averbakh2005computing} designs numerical methods to compute competitive ratios efficiently 
%A competitive ratio can be negative, and will be nonnegative if and only if $D(0) \le 0$ due to monotonicity of $D(\beta)$ from Theorem \ref{thm:slope-bounds}.

\subsection{Conservatism Control.} \label{sec:control}

As discussed earlier, overconservatism can be caused by criteria obsessed with worst-case scenarios while ignoring all other opportunities. 
The problem is exacerbated if ignored opportunities are from very likely scenarios, which may be furnished by experts. %consisting of most likely values
The ARM criterion offers a continuum of conservatism via $\beta$ by a family of robust policies $\pi^*_\beta$ from \eqref{defn:policy-form}, so that a choice of $\beta$ can be made to best catch such opportunities to improve performance and mitigate overconservatism.

Depending on the context of application, there may be different rationales on how to choose a suitable $\beta$.
Practitioners often desire robust policies with high expected reward, such as in the case of revenue management (see \citealt{vinod2021approach}). Therefore, an ideal $\beta$ should satisfy 
\begin{equation}\label{def:beta2max.expected}
\max_{\beta} E_{\zeta} \; r(\pi^*_\beta, \zeta).
\end{equation} 
%which provides a robust solution with least loss in expected reward.
Sometimes a robust policy with performance guarantee may be desirable for risk control, even when the distribution on scenarios is available.
In such case it is also appropriate to choose a policy $\pi^*_\beta$ according to \eqref{def:beta2max.expected}, which will be referred to as the robust optimal performance (ROP) policy.

Though the ROP policy is unavailable without a distribution, methods may be inspired to choose $\beta$ for conservatism control by adapting to likely opportunities near the mode of distribution, if unimodality is assumed. 
Experts may be asked to provide a typical scenario $\zeta^*$ made up of most likely values to represent likely opportunities.
If a policy performs well for $\zeta^*$, then it is likely to boost the expected reward in \eqref{def:beta2max.expected}, thus $\beta$ may be chosen by solving $\max_{\beta} r(\pi^*_\beta, \zeta^*)$, which can be difficult as $\pi^*_\beta$ may be very complicated.

An easier problem may result if opportunities are measured in terms of ex post optimal rewards.
Assuming unimodality for ex post optimal rewards, let $\hat{r}^*$ be the mode estimated by experts, which may serve as a {\em representative opportunity} $\ddot{r}$ of likely opportunities.
If $r(x,\zeta)$ is continuous and for some reason experts only provide the typical scenario $\zeta^*$, then a scenario $\zeta$ close to the typical scenario may satisfy $r^*(\zeta) \approx r^*(\zeta^*)$, thus $r^*(\zeta^*)$ may also serve as a representative $\ddot{r}$ just like $\hat{r}^*$. % , as $r^*(\zeta)$ must also be continuous
With a representative $\ddot{r}$ in place, consider a scenario family $U(\ddot{r})=\{\zeta\in\mathcal{U}:r^*(\zeta)=\ddot{r}\}$.
Without a distribution, the performance guarantee $\check{r}(\ddot{r}, \beta)$ for $U(\ddot{r})$ may serve as a proxy for the expected reward in \eqref{def:beta2max.expected}, so that a suitable $\beta$ can be chosen for the optimal performance guarantee (OPG):
\begin{equation}\label{def:beta-heuristic}
	\check{r}^*(\ddot{r}) = \max_{\beta} \check{r}(\ddot{r}, \beta) = \max_{\beta} \beta \ddot{r} - D(\beta).
\end{equation}
This framework of OPG can adapt the level of conservatism of the ARM criterion to likely opportunities represented by $\ddot{r}$, which roughly corresponds to scenarios in $U(\ddot{r}, \delta)=\{\zeta: r^*(\zeta)\in [\ddot{r}-\frac{\delta}{2}, \ddot{r}+\frac{\delta}{2}]\}$, where $\delta$ is a small diameter.
The OPG $\check{r}^*(\ddot{r})$ turns out to be the convex conjugate of $D(\beta)$, a construct quite amenable to analysis. %, usually nonnegative
Note that the OPG framework can be adopted by DRO if the ambiguity set is specified by confidence intervals of distribution parameters, where $\zeta^*$ would correspond to a typical distribution with parameters set to the nominal values or point estimates.

To analyze the properties of the OPG framework, let $r^*_- = \min_\zeta r^*(\zeta)$ and $r^*_+ = \max_\zeta r^*(\zeta)$, and let $\beta^*_-(r)=\min B^*(r)$ and $\beta^*_+(r)=\max B^*(r)$, where $B^*(r)=\{\beta: \beta r - D(\beta)=\check{r}^*(r) \}$.
When $r\in\{r^*_-, r^*_+\}$, an upper bound for $\check{r}(r, \beta)$ is 
\begin{equation}\label{def:ub4reward.guarantee}
\bar{r}(r) = \max_{\pi\in\Pi} \min_{\zeta \in U(r)} r(\pi,\zeta),
\end{equation}
which comes from 
\begin{eqnarray*}
D(\beta) &=& \min_{\pi\in\Pi}\max_{\zeta\in\mathcal{U}} \beta r^*(\zeta) - r(\pi,\zeta)\\
 &\ge& \min_{\pi\in\Pi}\max_{\zeta\in U(r)}\beta r^*(\zeta) - r(\pi,\zeta)\\
 &=&\beta r - \bar{r}(r).	
\end{eqnarray*}

\begin{thm}%[OPG Properties]
 \label{thm:best-guarantee}
If $r(x, \zeta)$ is bounded, the OPG $\check{r}^*(r)$ and the maximizers $\beta^*_-(r)$ and $\beta^*_+(r)$ in \eqref{def:beta-heuristic} have these properties:

i. For $r_1 > r_0$, there is $\beta^*_-(r_1) \ge \beta^*_+(r_0)$,
%i. Both the maximizer $\beta^*(r)$ and the related regret guarantee $D(\beta^*(r))$ increase in $r$ on $[0,\infty)$.
the OPG $\check{r}^*(r)$ is convex with $\check{r}^*(r_0) < \check{r}^*(r_1)$ if $\beta^*_+(r_0) > 0$ and 
$\check{r}^*(r_0) > \check{r}^*(r_1)$ if $\beta^*_-(r_1) < 0$, and
\[\begin{cases}
\beta^*_+(r)=-\infty,\; \check{r}^*(r)=+\infty & \text{if }  r<r^*_-\\
\beta^*_-(r)=-\infty,\; \check{r}^*(r) \le \bar{r}(r^*_-) & \text{if }  r=r^*_-\\
\beta^*_-(r)>-\infty,\; \beta^*_+(r)<+\infty & \text{if }  r\in(r^*_-,r^*_+)\\
\beta^*_+(r)=+\infty,\; \check{r}^*(r) \le \bar{r}(r^*_+) & \text{if } r=r^*_+\\ % \phantom{-}
\beta^*_-(r)=+\infty,\; \check{r}^*(r)=+\infty & \text{if }  r>r^*_+
\end{cases}\]
% with these special values
% \[\check{r}^*(r)=\begin{cases}
% \infty & \text{ if } r<r^*_-\\
% \displaystyle\max_\pi \min_{r^*(\zeta)=r^*_-} r(\pi,\zeta) & \text{ if } r=r^*_-\\
% \displaystyle\max_\pi \min_{r^*(\zeta)=r^*_+} r(\pi,\zeta) & \text{ if } r=r^*_+\\
% \infty & \text{ if } r>r^*_+
% \end{cases}\]

% $\check{r}^*(r)=\infty$ if $r<r^*_-$ or $r>r^*_+$, $\beta^*_+(r)=-\infty$ if $r<r^*_-$ and $\beta^*_-(r)=\infty$ if $r>r^*_+$, $\beta^*_-(r^*_-)=-\infty, \check{r}^*(r^*_-)=\max_\pi \min_{\zeta:r^*(\zeta)=r^*_-} r(\pi,\zeta)$ and $\beta^*_+(r^*_+)=\infty, \check{r}^*(r^*_-)=\max_\pi \min_{\zeta:r^*(\zeta)=r^*_+} r(\pi,\zeta)$, 

% and for any $r_0 < r < r_1$ with $\beta^*(r_0)<\beta^*(r)<\beta^*(r_1)$ there is  
%\[\frac{D(\beta^*(r))-D(\beta^*(r_0))}{\beta^*(r)-\beta^*(r_0)}\le r \le\frac{D(\beta^*(r_1))-D(\beta^*(r))}{\beta^*(r_1)-\beta^*(r)}.\]
% And for any $r_0 < r < r_1$: \frac{\check{r}^*(r)-\check{r}^*(r_0)}{r-r_0}\le \beta^*(r) \le\frac{\check{r_1}^*(r)-\check{r}^*(r)}{r_1-r}

ii. For $r\in[r^*_-, r^*_+]$, the absolute guarantee gap (AGG) of $G(r) = r - \check{r}^*(r)\ge 0$, is concave and strictly increases (decreases) in $r$ if $\beta^*_+(r)<1$ (if $\beta^*_-(r)>1$).
When $r\ge r^*_- > 0$, the relative guarantee gap (RGG) of $G(r)/r$ strictly increases (decreases) in $r$ if $\beta^*_+(r)<\gamma^*$ (if $\beta^*_-(r)>\gamma^*$).
% guess: the relative guarantee gap $G(r)/r$ peaks at the competitive ratio.
% if $D'(\beta)$ is continuously increasing, then (G/r)' = - D(\beta^*(r))/r^2.
\end{thm}
\begin{proof} These properties are proved as follows.
i. Let $\beta^*_i \in B^*(r_i)$ for $i=0,1$.
% Clearly there is $\check{r}^*(r_1) = \check{r}(r_1, \beta^*_1) \ge \check{r}(r_1, \beta^*_0) \ge \check{r}(r_0, \beta^*_0) = \check{r}^*(r_0)$.
Assume $\beta^*_1 < \beta^*_0$, optimality of $\beta^*_0$ gives $\check{r}(r_0, \beta^*_1) \le \check{r}(r_0, \beta^*_0)$  
\begin{eqnarray*}
	\Rightarrow &r_0 \beta^*_1 - D(\beta^*_1) &\le r_0 \beta^*_0 - D(\beta^*_0)\\
	\Rightarrow &D(\beta^*_0)- D(\beta^*_1) &\le r_0 (\beta^*_0 - \beta^*_1) \\
	\Rightarrow &D(\beta^*_0)- D(\beta^*_1) &< r_1 (\beta^*_0 - \beta^*_1) \\
	\Rightarrow &r_1 \beta^*_1 - D(\beta^*_1) &< r_1 \beta^*_0 - D(\beta^*_0),
\end{eqnarray*}
leading to $\check{r}(r_1, \beta^*_1) < \check{r}(r_1, \beta^*_0)$, a contradiction to the optimality of $\beta^*_1$, which proves $\beta^*_1 \ge \beta^*_0$, implying $\beta^*_-(r_1) \ge \beta^*_+(r_0)$ when $\beta^*_1 = \beta^*_-(r_1), \beta^*_0 = \beta^*_+(r_0)$.

Note that $\check{r}^*(r)$ is convex conjugate of $D(\beta)$, a pointwise maximum of a family of strictly increasing (decreasing) affine functions with slope $\beta^*_+(r_0) > 0$ ($\beta^*_-(r) < 0$) of $r$, hence convex and strictly increasing after $r_0$ (decreasing before $r_1$) with monotonicity of $\beta^*_-(r)$ and $\beta^*_+(r)$. % Fenchel conjugate

By Theorem \ref{thm:slope-bounds}, there is 
$r^*_- \le (D(\beta_1) - D(\beta_0))/(\beta_1-\beta_0) \le r^*_+$ 
for $\beta_1 > \beta_0$, which gives 
$D(\beta_0) \le D(\beta_1) - (\beta_1-\beta_0) r^*_-$, so that 
$\check{r}(r, \beta_0) \ge \check{r}(r, \beta_1) + (\beta_1-\beta_0) (r^*_- - r)$. 
Clearly, if $r < r^*_-$, there is $\lim_{\beta_0\downarrow -\infty}\check{r}(r,\beta_0)=\infty$, giving $\check{r}^*(r)=\infty$ and $\beta^*_+(r)=-\infty$. 
If $r = r^*_-$ then $\check{r}(r,\beta)$ increases as $\beta\downarrow-\infty$, thus $\beta^*_-(r)=-\infty$ and $\check{r}^*(r^*_-) = \lim_{\beta\downarrow-\infty} \beta r^*_- - D(\beta)$, which exists as it increases with an upper bound $\bar{r}(r^*_-)$.
% \begin{eqnarray*}
% \beta r^*_- - D(\beta) &=& \beta r^*_- - \min_{\pi\in\Pi}\max_{\zeta\in\mathcal{U}}\beta r^*(\zeta) - r(\pi,\zeta)\\
%   &\le& \beta r^*_-  - \min_{\pi\in\Pi}\max_{\zeta\in U(r^*_-)} \beta r^*(\zeta) - r(\pi,\zeta)\\
%   &=& \max_{\pi\in\Pi}\min_{\zeta\in U(r^*_-)} r(\pi,\zeta).
% \end{eqnarray*}
For the case of $r\in(r^*_-,r^*_+)$, let $\pi^*_\beta$ be the optimal policy, $\zeta^*_\beta$ the worst scenario, and $\zeta^*_- \in U(r^*_-)$. Clearly,  
\[D(\beta)= \beta r^*(\zeta^*_\beta) - r(\pi^*_\beta,\zeta^*_\beta) \ge \beta r^*(\zeta^*_-) - r(\pi^*_\beta,\zeta^*_-) \ge (\beta - 1) r^*_-.\]
Then $r(r, \beta) \le \beta (r - r^*_-) + r^*_-$, which goes to $-\infty$ as $\beta\downarrow-\infty$, thus $\beta^*_-(r)>-\infty$.
The results involving $r^*_+$ follow similarly.
% Let $\pi^*_\beta$ be the optimal policy and $\zeta^*_\beta$ the worst scenario, so that $D(\beta)= \beta r^*(\zeta^*_\beta) - r(\pi^*_\beta,\zeta^*)$, there is 
% \[\frac{\check{r}^*(r^*_-)}{\beta} 
% = r^*_- - \lim_{\beta\downarrow-\infty} \frac{D(\beta)}{\beta} 
% = r^*_- - \lim_{\beta\downarrow-\infty} r^*(\zeta^*_\beta),\] 

% \beta r - D(\beta_1) + (\beta_1-\beta_0) r^*_- =

% For the other claim, the optimality of $\check{r}^*(r')$ for any $r'\ge 0$ gives $\check{r}^*(r') \ge \beta^*(r) r' - D(\beta^*(r))$, or $D(\beta^*(r)) \ge \beta^*(r) r' -  \check{r}^*(r')$, thus
% \begin{equation*}
% 	D(\beta^*(r)) = \max_{r' \ge 0} \beta^*(r) r' -  \check{r}^*(r'),
% \end{equation*}
% a pointwise maximum of a family of increasing functions (slope $r'\ge 0$) of $\beta^*(r)$, hence increasing in $\beta^*(r)$, which in turn increases in $r$.

ii. Rewrite $G(r)$ into $G(r) = \min_{\beta\ge 0} D(\beta) + (1-\beta) r$, a pointwise minimum of a family of strictly increasing (decreasing) affine functions in $r$ when $\beta^*_+(r)<1$ ($\beta^*_+(r)>1$).
Rewrite $G(r)/r$ for $r>0$ into $G(r)/r = \min_{\beta\ge 0} D(\beta)/r - \beta +1$, a pointwise minimum of a family of increasing (decreasing) functions in $r$ when $D(\beta)<0$ ($D(\beta)>0$), as $D(\beta)$ strictly increases with $r^*_->0$ by Theorem \ref{thm:slope-bounds}. 
By Theorem \ref{thm:cr-equivalence}, when $\beta^*_+(r)<\gamma^*$ ($\beta^*_-(r)>\gamma^*$), there is $D(\beta^*_+(r))<0$ ($D(\beta^*_-(r))>0$), and $G(r)/r$ should strictly increase (decrease) in the neighborhood of $r$.
\end{proof}

Theorem \ref{thm:best-guarantee} reveals some nice properties of the ARM criterion and the OPG framework. 
It shows that to adapt to opportunities concentrated around a bigger $r$, the ARM criterion will employ a bigger $\beta^*$, and adjust its level of conservatism accordingly.
As there is always $r \in [r^*_-, r^*_+]$, the OPG $\check{r}^*(r)$ is a bounded convex function, while the sign of $\beta^*_-(r)$ and $\beta^*_+(r)$ can help determine if $\check{r}^*(r)$ increases or decreases.
%Though the $\beta$-adjusted regret guarantee $D(\beta^*(r))$ also increases with $r$, 
The AGG and RGG are indicators of efficiency to capture opportunities, where RGG is more appropriate for relative losses.
High efficiency (with small AGG or RGG) when $r$ is near the ends of $[r^*_-, r^*_+]$ may be interpreted as opportunities are easier to capture when they are cornered to either ends.
If $D(\beta)$ is continuous and convex, then $D(\beta)$ and $\check{r}^*(r)$ are mutually convex conjugate.
%The properties of $D(\beta)$ are further studied for better use of the heuristic.
% These results also suggest that the conservatism may change in a ``smooth'' manner with $\beta$, as well as the optimal policy.
	%\begin{cor}
	%If \(r^*(\zeta) > 0\) for all $\zeta \in \mathcal{U}$, there is a unique solution to $D(\beta)=0$ on the interval $[-1, 1]$.
	%\end{cor}
	%
	%\begin{proof}
	%	It follows from the continuity and monotonicity of $D(\beta)$ and the simple fact of $D(-1) \leq 0 \leq D(1)$.
	%\end{proof}
%\subsection{Convexity.}
%The convexity of $D(\beta)$ requires certain conditions to hold. 
Note that $\beta r^*(\zeta) - r(\pi,\zeta)$ is linear in $\beta$, thus 
\(  F(\beta; \pi) = \max_{\zeta\in\mathcal{U}} 
	\beta r^*(\zeta) - r(\pi,\zeta)
\) 
is convex in $\beta$ for a given policy $\pi$. 
But $D(\beta) = \min_{\pi\in\Pi} F(\beta; \pi)$ is not necessarily convex in $\beta$, and extra conditions are required to ensure its convexity.

\begin{thm}%[$\beta$-Convexity]
\label{thm:beta-convexity}
With consistent convexity, the best regret guarantee $D(\beta)$ is convex, and strictly so if $r(x,\zeta)$ is strictly concave in $x$ for all $\zeta\in\mathcal{U}$.
\end{thm} 
\begin{proof} 
Let $\ddot{\pi}^*_i\in \ddot{\Pi}$ be an optimal pruned policy for $\beta_i, i=1,2$, 
and let $\ddot{\pi}' = \lambda \ddot{\pi}_1 + (1-\lambda)\ddot{\pi}_2$ for any $\lambda \in (0,1)$.
By the concavity of $r(x,\zeta)$, there is $r(\ddot{\pi}',\zeta) \geq 
	\lambda r(\ddot{\pi}^*_1,\zeta) + (1-\lambda) r(\ddot{\pi}^*_2,\zeta)$
for all $\zeta\in\mathcal{U}$.
Let $\beta = \lambda \beta_1 + (1-\lambda) \beta_2$ and proceed as follows:
		\begin{eqnarray*}
			D(\beta) &=& \min_{\ddot{\pi}\in \ddot{\Pi}} \max_{\zeta\in\mathcal{U}} 
			\beta r^*(\zeta) - r(\ddot{\pi},\zeta)\\
			&\leq& \max_{\zeta\in\mathcal{U}} \beta r^*(\zeta) - r(\ddot{\pi}',\zeta) \\
			&\leq& \max_{\zeta\in\mathcal{U}} \beta r^*(\zeta) - \left(\lambda r(\ddot{\pi}^*_1,\zeta) + (1-\lambda) r(\ddot{\pi}^*_2,\zeta) \right)\\
			&\leq& \lambda \left(\max_{\zeta\in\mathcal{U}} \beta_1 r^*(\zeta) - r(\ddot{\pi}^*_1,\zeta)\right) +\\
			&& (1-\lambda) \left( \max_{\zeta\in\mathcal{U}} \beta_2 r^*(\zeta) - r(\ddot{\pi}^*_2,\zeta)\right)\\
			&=&\lambda D(\beta_1) + (1-\lambda) D(\beta_2).
		\end{eqnarray*}
Similarly, strict concavity comes by $r(\ddot{\pi}',\zeta) > 
	\lambda r(\ddot{\pi}^*_1,\zeta) + (1-\lambda) r(\ddot{\pi}^*_2,\zeta)$. 
\end{proof}

%As the benchmark is linearly scaled up by $\beta$, the convexity of $D(\beta)$ indicates that as $\beta$ increases, the corresponding optimal solution is less and less speedy to catch up with the scaling up of the benchmark. 

%Theorem \ref{thm:beta-convexity} ensures that local optimality of \eqref{def:beta-heuristic} is global, and it is unique with strict convexity.
If $D(\beta)$ is strictly convex and differentiable with a continuous and strictly increasing $D'(\beta)$, then $\beta^*(r)$ for \eqref{def:beta-heuristic} is simply the inverse of $D'(\beta)$ from the first order condition
\begin{equation}\label{eqn:beta-derivative}
\frac{\partial \check{r}(r,\beta)}{\partial \beta} = r - D'(\beta) = 0.
\end{equation}
%Let $\beta^*$ solves \eqref{eqn:beta-derivative}. 

%which can be interpreted as the ARM criterion gets less conservative for more opportunities in $\zeta^*$.

This section illustrates how the ARM criterion adapts the level of conservatism to opportunities, while maintaining all advantages of RO.
Only a representative opportunity is required by the OPG framework as additional information, and a new approach to competitive ratio is also established.
It is time to see these theoretical results in action.

\section{One-way Trading.}\label{sec:1-way}

% TODO: cite a survey paper on 1-way trading
%It has been studied extensively using various tools like stochastic optimization, dynamic programming, and reinforcement learning. 
%The first work under RO with competitive ratio is \cite{el2001optimal}, which aroused much interest. 
%\cite{wang2016competitive} follows up the work and employs the absolute regret criterion with closed-form solutions. 
%while recovering the main results of both \cite{el2001optimal} and \cite{wang2016competitive} by setting $\beta$ to a proper value.

In this section the ARM criterion is applied to the one-way trading problem to demonstrate its properties and potential, such as the effectiveness of conservatism control and the new approach to competitive ratio analysis. 
The one-way trading problem has been richly studied with both competitive ratio (\citealt{el2001optimal}) and absolute regret (\citealt{wang2016competitive}), which are ideal for comparison.
Closed-form analytical solutions are derived with the ARM criterion, yielding the result of \cite{wang2016competitive} as a special case with $\beta=1$. %while the derivation process is not much more complex than theirs.
The competitive ratio is directly found via the new approach, in contrast to the intricate derivation in \cite{el2001optimal} that heavily depends on acute intuition and deep insights of great talents. % but are unnecessary with the new approach
Finally, numerical simulations are conducted to verify that the ARM criterion can indeed offer smooth control of conservatism, and the OPG framework can determine a suitable level of conservatism with improved performances. 
% These are helpful to demonstrate the properties and potential of the ARM criterion.
	
\subsection{Problem Formulation.}
Consider selling a certain amount of divisible goods (such as gasoline) in  periods $t = 1,\cdots,T$, while the price fluctuates in the range of $[m, M]$. 
In each period, a single price $p_t \in [m, M]$ is revealed first, then the trader sells as a price-taker at $p_t$ an amount $x_t \ge 0$ out of the remaining stock, without knowing the future prices. %commit to 
In the last period $T$, the trader must sell out whatever remains. 
The goal is to maximize the total sales revenue.
	
This is a multi-stage problem, in which the stages naturally coincide with periods.
A scenario $\zeta$ is simply the prices $p=(p_1,\cdots, p_T)$ revealed over time, with $\zeta_t = p_t$. 
There is $\mathcal{U}=[m,M]^T$ and $\mathcal{U}_t(\zeta_{1:t}) = [m,M]$, with prices independent of each other.
Without loss of generality, the total amount to sell is 1 unit, and the action is $x = (x_1, \cdots, x_T)$ with $X = \{x \ge 0: \sum_{t=1}^T x_t = 1\}$. 
For $t<T$ there is $X_t(h_t) = [0, q_t]$ where $q_t = 1-\sum_{s=1}^{t-1} x_s$ is the remaining stock to sell given $h_t$, but in the last period $X_T(h_T) = [q_T, q_T]$.
The reward is accrued over time, so let $r_t = \sum_{s=1}^{t-1} p_s x_s$ for the rewards accrued over $h_t$, and the reward in the end is $r(x, p) = r_{T+1}$.
Let $\hat{p}_t = \max\{ p_s: s=1,\cdots, t-1\}$ denote the highest price in $h_t$, and $r^*(p) = \max \{r(x,p): \sum_{t=1}^T x_t = 1 \} = \hat{p}_{T+1}$ the ex post optimal. 
In the end period (\ref{defn:D[T](h)}) becomes
	\begin{equation} \label{1way:D[T](h)}
		D_T(h_{T+1}; \beta) = \beta \hat{p}_{T+1}  - r_{T+1}.
	\end{equation}
	
By tradition, in a stage $t$ the price $p_t$ is revealed first, then an action $x_t$ is taken, which differs from the standardized formulation in \eqref{defn:Dt(h)}: %from the standardized in (\ref{defn:D[t](h)}): 
\begin{eqnarray}\label{1way:D[t](h)}
		D_{t-1}(h_t; \beta)
		&=& \max_{p_t \in [m,M]} \min_{x_t \in X_t(h_t)}  D_t(h_{t+1};\beta).
\end{eqnarray}
As a dummy decision can be inserted in the beginning to standardize it, this difference is superficial, and all results in Section \ref{sec:prop} remain valid.

\subsection{Analytic Solution.}
The analysis starts from the last period $T$ and works backwards.
It is enough to consider $\beta\ge 0$ for comparison with related work, while keeping it simple. 
In the last period there is $x_T=q_T$, and (\ref{1way:D[t](h)}) becomes
\[
	D_{T-1}(h_T; \beta) = \max_{p_T\in[m,M]} 
	\beta \max(\hat{p}_{T}, p_T) - (r_{T} + p_T q_T),
\]
which is convex in $p_T$, and the maximizer is either $p_T=m$ or $p_T=M$, thus
\begin{eqnarray*}
	D_{T-1}(h_T; \beta) &=& \max(\beta \hat{p}_{T} -  R_{T}, \beta M - (r_{T} + M q_T))\\
	&=& \max(\beta \hat{p}_{T}, \beta M - (M-m) q_T) -  R_{T}\\
	&=& \beta \max(\hat{p}_{T}, P_1(q_T)) - R_{T},
\end{eqnarray*}
where $R_{t} = r_{t} + m q_{t}$ for $t = 1,\cdots,T$ is the lower bound on $r_{T+1}$ given $h_t$, and $P_1(y)$ is an auxiliary quantity-to-price function defined as
\[ P_j(q) = (M-m)\left(1-{q\over\beta j}\right)_+^{j}+m, j=1, 2, \cdots, \]
with $y_+^{j} = \max^j(0,y)$ for the positive part of $y$ raised to the $j^{th}$ power. 
Let $P^-_j(y) = q$ be the inverse of $y = P_j(q)$ for $q\in [0,\beta j]$, with the trivial case of $\beta=0$ found to be $P^-_j(y) \equiv 0$ in limit as $\beta\downarrow 0$.
Continue on with (\ref{1way:D[t](h)}) for $t=T-1, \cdots, 1$ by backward induction, it is solved analytically:

\begin{thm} \label{thm:1way-AR}
The minimal worst-case regret for the one-way trading problem in period $t$ given history $h_t$ for $t=1, 2, \cdots, T$ is
\begin{equation}\label{eq:Dt}
		D_{t-1}(h_t; \beta) 
		= \beta\max(\hat{p}_{t}, P_{1+T-t}(q_t))-R_{t},
\end{equation}
and the optimal trading policy is $\pi^*_t(h_t, p_t) = q_t - q_{t+1}^*$, where $q_{T+1}^*=0$ and 
\begin{equation}  \label{eqn:1way-q*-solved}
		q_{t+1}^* = \min(q_{t}, P^-_{T-t}(\hat{p}_{t+1})), ~ t=1, \cdots, T-1.
\end{equation}
\end{thm}

The result of \cite{wang2016competitive} is a special case of Theorem \ref{thm:1way-AR} with $\beta=1$, and the proof for Theorem \ref{thm:1way-AR} requires more general treatments, but the complexity of analysis is similar, which is presented in \ref{app:1way-AR}.
The next three corollaries are readily derived from Theorem \ref{thm:1way-AR}.
	
	%Corollary: the worst price and optimal quantity are
	%\begin{equation}
	%q^*_t = 
	%p^*_{t+1} =
	%\end{equation}
	
\begin{cor} %[Convex Regret Guarantee]
\label{cor:1way-AR}
The minimal worst-case regret $D(\beta)$ for the one-way trading problem is a convex function of $\beta$:
	\begin{equation}\label{eq:DT}
		D(\beta) = \beta (M-m)\left(1-{1\over \beta T}\right)_+^{T} - (1-\beta) m,
	\end{equation}
\end{cor}
\begin{proof} 
In the first period, there is $q_1=1, r_{1}=0, \hat{p}_{1}=m$. Use these in (\ref{eq:Dt}) and simplify to have the result. The convexity of $D(\beta)$ is a consequence of the consistent convexity of the one-way trading problem and Theorem \ref{thm:beta-convexity}.
\end{proof}

Based on Theorem \ref{thm:1way-AR}, the competitive ratio can be directly derived from the new approach of competitive analysis, which is much simpler than the truly ingenious but intricate analysis of \cite{el2001optimal}. %Such simplification is impossible with the analysis and results in \cite{wang2016competitive}.

	%Corollary: there are $2^n$ equilibrium in the adversary-trader game.
\begin{cor}%[Competitive Ratio]
The competitive ratio defined in (\ref{defn:cratio}) for the one-way trading problem is the unique root of $D(\beta)=0$ for $D(\beta)$ defined in (\ref{eq:DT}).
\end{cor}
\begin{proof} 
As $r^*(\zeta) \geq m > 0$, it follows from Theorem \ref{thm:cr-equivalence}.
\end{proof} 

This result agrees perfectly with \cite{el2001optimal}, except that they define competitive ratio as the inverse of the definition employed here. 
Their analysis is much more involved and heavily relies on insights of the worst case price paths, which is not required here. %which can be derived in a way similar to what is done in \cite{wang2016competitive}.

\begin{cor} % [Increasing Aggressiveness]
\label{cor:increasing-aggressiveness}
As $\beta$ increases, the optimal trading policy gets more optimistic and aggressive by trading less for now and more for the future to take on more risks under the same $h_t$ and $p_t$, which means that for $\beta_1 > \beta_2 > 0$ and $\forall h_t\in H_t, \forall p_t\in [m,M]$, there is \( \pi^*_t(h_t,p_t; \beta_1) \le \pi^*_t(h_t,p_t; \beta_2) .\)
\end{cor}
\begin{proof} 
Consider the quantity reserved for future $q^*_{t+1}$ in \eqref{eqn:1way-q*-solved} and note that 
\[
		P^-_{T-t}(p) = \beta (T-t) \left( 1 - \sqrt[T-t]{p-m\over M-m} \right)
\]
increases in $\beta$, therefore $q^*_{t+1}$ increases as $\beta$ increases.
\end{proof} 

Corollary \ref{cor:increasing-aggressiveness}  analytically illustrates the continuous moderation of conservatism by $\beta$: as $\beta$ increases, the optimal policy becomes less conservative.
% by trading less for now and more for the future under the same condition of $(h_t,p_t)$. 

\subsection{Numerical Study.}
The effect of conservatism control is further studied numerically by counterfactual simulations on the one-way trading problem.
The setup is similar to that of \cite{wang2016competitive}, the main difference is in price distribution: instead of uniform distribution, Beta distribution is employed to allow for variations of most likely opportunities.
The instance has $T=100$ periods, and the price range is $[m,M]=[60, 120]$, which is also the range $[r^*_-, r^*_+]$.
The prices for all periods are independent and identically distributed (IID) with a Beta(2, 2) distribution scaled over $[m, M]$, which is symmetric and unimodal at the center.
Suppose the experts provide an accurate estimation of $\hat{r}^*=117.5$, so that
the OPG with the representative $\ddot{r}=\hat{r}^*$ is obtained at $\beta=3.14$ by numerically solving \eqref{eqn:beta-derivative}, which will be labeled ``MLO'' for most likely opportunity.

% $\hat{r}^*=2.897$, calcuated from $r^*(\zeta)=\hat{p}_{T+1} = \max_{i=1}^T p_i$
% giving a mean of $0.5+a/(a+b)=1.1$, a variance of $ab/((a+b)^2(a+b+1))=0.04$ (or a standard deviation of 0.2), and a mode of $0.5+(a-1)/(a+b-2)=1.167$. %If a,b>1, (or one of them =1), 
% for n iid with CDF F or pdf f, the mode of the max solves: f' = - (n-1)f^2/F 
%but the ARM formulation only knows the price bounds $m=1, M=2$. 

\begin{figure}
	    \centering
	    % simfun(5,3,2,4, 1,2, 10000)
	    \includegraphics[scale=0.75]{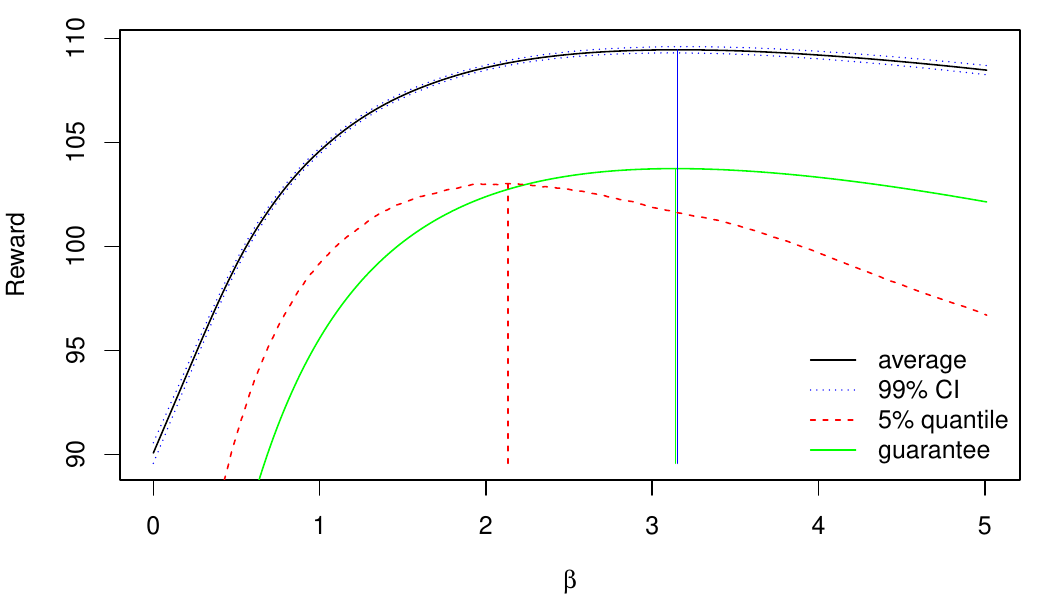}
	    \caption{Performance guarantee and average reward with 99\% CI.}
	    \label{fig:my_fig}
	\end{figure}

The simulation is conducted with the technique of common random numbers for variance reduction to better compare the average rewards of different policies.
A single run of the simulation would generate $T$ IID random prices first, and then feed them as common input to all the policies, with a total reward produced by each policy in the end.
Such a single run is repeated $N=5,000$ times, at the end of which the 5\% quantile and the average reward with its 99\% confidence interval (CI) are computed for each policy. %$\pi^*_\beta$. 
After the optimal policies $\pi^*_\beta$ in \eqref{eqn:1way-q*-solved} for $\beta \in [0,5]$ are evaluated, the ROP $\beta=3.16$ with a reward of $109.5$ is found from the average reward. 
These are plotted in Fig. \ref{fig:my_fig}, together with the performance guarantee $\check{r}(\hat{r}^*, \beta)$. 

Some interesting observations can be made on Fig. \ref{fig:my_fig}.
First, the performance guarantee $\check{r}(\hat{r}^*, \beta)$ is a fairly good proxy for the average reward, as their maximizers (the ROP $\beta=3.16$ and the MLO $\beta=3.14$) are very close.
%though there is a fairly big gap between the two curves.
Second, the OPG framework tends to reduce the downside risks, as observable by the 5\% quantile:
The quantile at the MLO $\beta$ is $101.7$, which is very close to the maximal quantile of $103.1$ (at $\beta=2.15$).
Thus, the MLO policy derived from the OPG framework in the experiment indeed improves performance, reduces conservatism, and maintains a relatively low risk.

% the heuristic works quite well, with 
Besides the optimal policies for $\beta \in [0,5]$, the simulation also includes the policy that maximizes expected reward (shorthanded as `max expected', see Appendix B of \citealt{wang2022robust} for technical details) and the ex post optimal policy for comparison. 
Seven policies are listed in Table \ref{table:policies},
%the maximin ($\beta=0$), the relative regret ($\beta=\gamma^*=0.78$), the absolute regret ($\beta=1$), the MLO ($\beta=3.14$), the ROP ($\beta=3.16$), the max expected, and the ex post optimal.
and the absolute and relative gaps are calculated by benchmarking against the `max expected' policy.
Note that the MLO policy improves the average reward by roughly 17\% over the maximin policy, 6\% over the relative regret policy, and 4\% over the absolute regret policy.
This kind of improvement can make a big difference in applications like the airline revenue management where the profit margin is usually razor-thin (\citealt{vinod2021approach}).
%Note that the heuristic makes more significant improvement over the other policies if the heuristic $\beta$ is further away from theirs, which helps choose the values of $a=3.5, b=1.5$ for Beta(a,b) to show off the potential of the ARM criterion and the heuristic.

%If there are more periods to sell, it pays to be more aggressive, and one would expect that the transition to the next phase will come later.
%on $\beta$, which is indeed the case as in Fig. \ref{fig:my_fig2} with $T=15$. 

%	\begin{figure}
%	    \centering
%	    % makeplot(10000, 15, (1:300)/100)
%	    \includegraphics[scale=0.75]{arc-sim-T15.eps}
%	    \caption{Average revenue and standard deviation of revenue for $T=15$}
%	    \label{fig:my_fig2}
%	\end{figure}

\begin{table}%[h!] 
	\centering
\begin{tabular}{c || r@{$\pm$}l | c | c } \hline
Policy ($\beta$) & Average & 99\% CI & Gap & Gap\% \\ \hline
maximin reward (0.00) & 90.267 & 0.348 & 25.469 & 22.0\% \\ \hline
relative regret (0.78) & 102.703 & 0.123 & 13.033 & 11.3\% \\ \hline
absolute regret (1.00) & 104.593 & 0.104 & 11.144 & 9.6\% \\ \hline
MLO (3.14) & 109.466 & 0.106 & 6.270 & 5.4\% \\ \hline
ROP (3.16) & 109.466 & 0.107 & 6.270 & 5.4\% \\ \hline
max expected & 115.737 & 0.073 & 0.000 & 0.0\% \\ \hline
ex post optimal & 116.859 & 0.043 & -1.122 & -1.0\% \\ \hline
\end{tabular}
\caption{Benchmark the ARM policies of various $\beta$ values.}
\label{table:policies}
\end{table}

% It seems there is a sweet spot in between extreme conservatism and aggressiveness.
% The extremely aggressive case of $\beta\uparrow\infty$ (absent in Fig. \ref{fig:my_fig}) almost only sells in the last period by \eqref{eqn:1way-q*-solved}, while the extremely conservative case of $\beta=0$ almost only sells in the first period, giving them identical average reward and standard deviation.
% Thus both extremes give poor performance with low expected reward and high overall risk. %with an optimal 
% As the standard deviation is unimodal with a minimal value of $0.179$ at $\beta=1.45$, the ARM criterion may achieve both higher rewards and lower risks sometimes with a proper $\beta$ value.
% $\beta$ depending on the context. 

%This can be gleaned from the following observations. 
% An overly conservative robust policy by the maximin reward criterion (i.e. $\beta=0$) can suffer both lower performance and higher overall risk. 
% On the other hand, an overly aggressive robust policy for a big $\beta$ may also hurt the performance as it boldly exposes to future price risks by reserving too much inventory, ending up selling it in the last period at any price. 

\begin{figure}
	    \centering
	    % broader.perspective(5, 1, 3)
	    \includegraphics[scale=0.75]{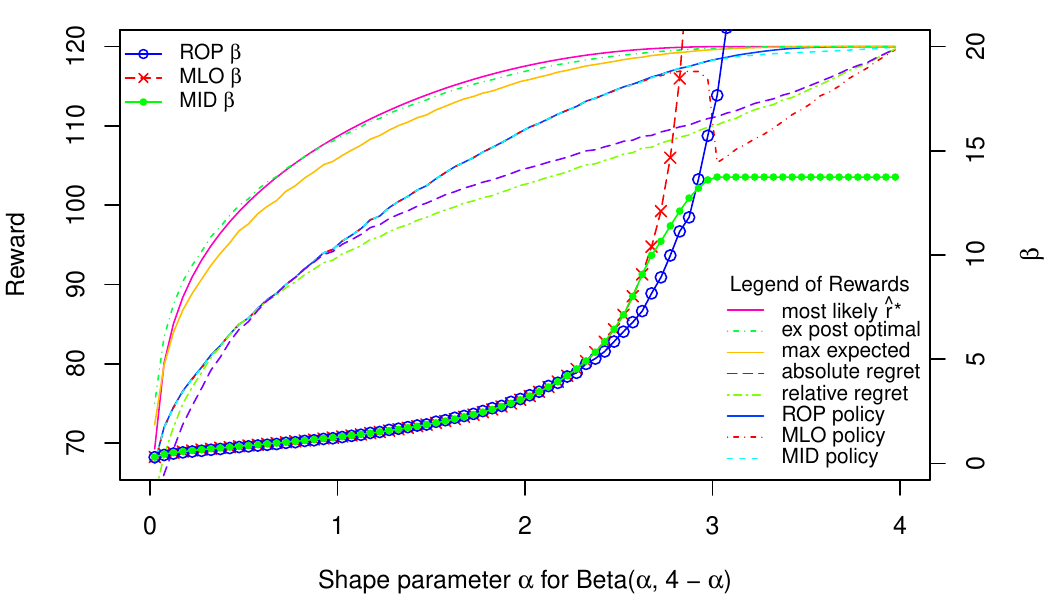} %% SAVE
	    \caption{Policies under a continuous shape-shifting scheme for the Beta distribution to simulate environments with increasingly better opportunities.}
	    \label{fig:broader}
	\end{figure}

To observe the effect of conservatism control in a broader perspective, the same experiment is conducted with a continuous shape-shifting scheme in Beta$(\alpha, 4-\alpha)$ for $\alpha\in (0, 4)$, with the end points excluded as they give degenerated distributions. 
Note that as $\alpha$ increases, the mode of Beta$(\alpha, 4-\alpha)$ moves toward $M$, and the opportunities concentrate more and more on the higher end of $[r^*_-, r^*_+]$ with the most likely opportunity $\hat{r}^* \uparrow r^*_+$.
Fig. \ref{fig:broader} presents the related rewards and $\beta$ values. 
Note that as $\alpha$ increases, both the MLO and ROP $\beta$ increase as expected,  giving more aggressive policies for more promising opportunities. 
%But the performance of MLO is not always ideal, so a detailed analysis is in order.
The average reward of MLO and ROP policies are hardly distinguishable for $\alpha \in (0, 2.8]$, which claims great successes for the MLO policies and the OPG framework in effective conservatism control with performance improvements. 

% A bigger $\alpha$ increases the probability for higher prices, there should be higher average reward and most likely opportunities $\hat{r}^*$, by which a bigger $\beta^*$ for a more aggressive or less conservative heuristic policy can be found from \eqref{eqn:beta-derivative}.
% and the heuristic $\beta^*$ gets out of simulation range if $b\downarrow 1$.
% Both the ROP and OPG $\beta$ steadily increases as $\alpha$ gets bigger in agreement to Theorem \ref{thm:best-guarantee}, illustrating smooth conservatism control to catch increasingly better opportunities.
% The MLO policy dominates the absolute and relative regret policy and performs very well relative to the ROP policy: The average rewards for the ROP and MLO policy are almost identical for $\alpha<2.7$, then the MLO policy starts to slightly fall behind, with a maximum difference of $0.68$ at $\alpha=2.9$ (as labeled out by a red diamond). 
%J-shaped with a left tail, convex
%negatively skewed, mode = 1
% In a beta distribution, which is a continuous probability distribution, the probability of the random variable being exactly equal to any single specific value, including 1, is 0.

Then a noticeable performance gap appears between the MLO and ROP policies, which widens as $\alpha$ increases on $(2.8, 3)$, with an abrupt MLO performance drop as $\alpha$ approaches 3. 
This gap then gradually narrows as $\alpha$ increases on $[3,4)$. 
This phenomenon can be explained by the discrepancy between the ROP and MLO $\beta$, which is significant after $\alpha=2.8$, with the MLO $\beta \uparrow \infty$ and $\hat{r}^*\uparrow r^*_+$ as $\alpha\uparrow 3$, while the ROP $\beta\approx 15$ at $\alpha=3$.
For $\alpha\in [3,4)$, the MLO $\beta=\infty$ with $\hat{r}^* = r^*_+$, while the ROP $\beta\uparrow\infty$ as $\alpha\uparrow 4$, converging to the MLO $\beta=\infty$.
Note that the policy $\pi^*_\beta$ for $\beta=\infty$ almost only trades in the last period, so that its expected reward is simply $60+15\alpha$, which is in agreement with simulated results in Fig. \ref{fig:broader}. 
% and the OPG $\beta$ then increases faster than and rapidly overshoots the ROP $\beta$, 
% Such phenomena in the rewards come from the $\beta$ values: The ROP and the OPG $\beta$ are relatively close in value for $\alpha\in (0, 2.8]$ with $\hat{r}^*\le 119.9$, then the difference rapidly widens as $\alpha\uparrow 3$ and the OPG $\beta\uparrow\infty$ and remains $\infty$ thereafter due to $\hat{r}^*=120 = M = r^*_+$ for $\alpha\in [3,4)$, while the ROP $\beta$ is about $15$ for $\alpha=3$ and approaches $\infty$ as $\alpha\uparrow 4$.

The inferior performance of the MLO policy with $\beta\uparrow\infty$ needs more analysis, though the situation of $\hat{r}^* \uparrow r^*_+$ is arguably rare in practice, and $\beta\uparrow\infty$ is not always the case in Theorem \ref{thm:best-guarantee}. %Such discrepancy is probably due to that the density of Beta$(\alpha, 4-\alpha)$ approaches $\infty$ at $M$ as $\alpha\uparrow 3$, and 
Note that at $\alpha=2.8$ the most likely $\hat{r}^* = 119.9$ is already very close to $r^*_+=120$, which left little room on the right of $\hat{r}^*$ to have it sit in the middle of likely opportunities to be representative of them.
% As $a\uparrow 4$, $\hat{r}^*\uparrow M$ and $M$ is the divergence point of probability density. 
%With $\hat{r}^*=M$ there is $\beta^*=\infty$ by Theorem \ref{thm:best-guarantee}, which explains the rapid increase in $\beta^*$.
%Note that such a disaster of divergence is peculiar to the one-way trading problem due to the ex post optimal being the maximum price over all periods. Nevertheless, it raises a red flag that feeding the heuristic with a $\hat{r}^*$ that is not so representative to proxy expected reward can negatively impact the performance.
A quick fix may consider the set of representative scenarios $U(\hat{r}^*, \delta)$, and the OPG $\beta$ can be determined from \eqref{eqn:beta-derivative} with the representative $\ddot{r}$ set to the midpoint of the supporting interval of $U(\hat{r}^*, \delta)$: $\ddot{r}=\frac{1}{2}(\max(r^*_-, \hat{r}^*-\frac{\delta}{2})+\min(r^*_+, \hat{r}^*+\frac{\delta}{2}))$.
In Fig. \ref{fig:broader} this fix with $\delta= 1\% (r^*_+ - r^*_-)$ is labeled ``MID'', whose average reward is comparable to that of the ROP policy throughout, with a maximal gap less than $0.6$.
Thus, the key for the OPG framework to provide effective conservatism control with improved performance is to have a representative $\ddot{r}$ surrounded by likely opportunities.
The value for the MID $\beta$ increases in $\alpha$ and then plateaus for $\alpha\in[3,4)$, which suggests simply setting an upper limit on the MLO $\beta$ heuristically as another possible fix.
There may still be other fixes when the MLO $\beta\uparrow\infty$, but the best practice probably depends on the application, which leaves an interesting topic for future research.

% As illustrated in Figure \ref{fig:my_fig}, the fine-tuned $\beta$ may well go beyond $\beta=1$ for the minimax absolute regret criterion.
% value depends on the problem context, which 
% A robust counterpart for standard deviation could be the reward spread of a policy, or the gap between its maximum and minimum reward obtained in all scenarios.
% If a $\beta$ can be identified whose corresponding optimal policy has the minimum reward spread, it might serve as a proxy for the the standard deviation minimizer in Fig. \ref{fig:my_fig}.

\section{Conclusion.}\label{sec:conclusion}

The ARM criterion proposed in this paper can provide fine control of conservatism by choosing a CCP ($\beta$) to catch likely opportunities, while maintaining all the major advantages of RO. 
% It is based on the idea of adjustable regret that benchmarks the performance of a solution against the ex post optimal performance multiplied by the control parameter.
This criterion minimizes the $\beta$-adjusted regret guarantee, from which a performance guarantee for any scenario can be obtained.
Convexity is preserved even for multi-stage problems, which is great for computational tractability and theoretical analysis.
Distributions are not required, experts only need to specify an opportunity representative of likely opportunities for conservatism control.
The performance guarantee for the representative opportunity is optimized in the OPG framework to calibrate $\beta$ to  catch opportunities and improve performance for conservatism control.
% The method is based on the mechanism that as $\beta$ increases, the ARM criterion will recommend solutions with better performance guarantees for scenarios of more opportunities. 
%is  the ARM criterion chooses a solution that is likely to ``mimic" the benchmark's  behavior, so that the conservatism of the recommended solution is moderated. 
% The control parameter happens to interpolate among well-known robust criteria of distinct levels of conservatism.
% helps choose a more aggressive solution with a larger control parameter that produces a more aggressive benchmark.
Various theoretical properties of the ARM criterion are studied, such as continuity, monotonicity, and convexity, to facilitate problem analysis, find closed-form solutions, or design better numerical algorithms.
The theoretical study also leads to a new approach for competitive ratio analysis, which may be simpler than the traditional approach, as in the analysis of the one-way trading problem.
Though reward maximization is considered here, parallel results can be readily derived for cost minimization.
% This new approach is successfully applied to the one-way trading problem to obtain the competitive ratio in a much simpler way.
% Though it is not clear what kind of problems are easier to tackle with this new approach to competitive ratio analysis, it seems from limited observations with the one-way trading problem that problems with additive objectives are more likely to benefit from this new approach, while problems with multiplicative objectives are less likely.

% The tractability is studied for two-stage linear problems with the ARM criterion. 
% Two particular situations are studied: the right-hand side uncertainty and the objective uncertainty.
% Equivalent reformulations into TSLRO/FR problems make it possible to take advantage of the tractable solution methods developed recently to solve them for practical applications. 

The ARM criterion is applied to the robust one-way trading problem to demonstrate its potential. Closed-form solution is obtained, from which the competitive ratios is easily derived by the new approach.
Analysis of the closed-form solution shows that the optimal policy gets more aggressive as $\beta$ increases.
Numerical experiments on one-way trading problem are designed to illustrate effective conservatism control, with significant performance improvements over other commonly used criteria.
%Either extreme conservatism or aggressiveness may suffer lower average reward and higher overall risk, and a proper $\beta$ from the heuristic or its variant can better catch the opportunities and deliver improved average reward. %, which depends on the problem context

This study of the ARM criterion may serve as a starting point for future researches. 
%Surely there are many directions for future study. 
% A rigorous analysis on how $\beta$ moderates conservatism is theoretically interesting and challenging at the same time. 
% and it can be used in many other applications as well, especially in situations that needs fine-tuned control over the conservatism of solution.
Applying it to practical problems in various areas is the thrust for further theoretical development. 
Note that it can also be applied when distributions can be estimated, but a performance guarantee is highly desired.
%It is especially of practical interest in revenue management
Innovative methods or frameworks may be developed to determine a suitable $\beta$ for conservatism control, which may lead to substantial contributions. 
%Researches on linear problems with the ARM criterion may be fruitful, as progresses are made in this regard on the absolute and relative regret criterion by \cite{poursoltani2021adjustable}.
Finally, it can be practically and theoretically fruitful to apply the ARM criterion with DRO.
% It is possible to adapt the ARM criterion and the heuristic method to DRO.

%Finally, the unimodality of variance in Fig. \ref{fig:my_fig} could also be a future research topic for applications interested in variance reduction.
%Finally, the unimodality of the standard deviation with regard to $\beta$ in Fig. \ref{fig:my_fig} suggests a possibility to minimize reward variance by choosing $\beta$, which might be a valuable direction of future study.

\appendix

\section{Proof of Lemma \ref{lem:optimality}}
\label{app:optimality}
\begin{proof}
A policy $\pi^* \in \Pi$ satisfying \eqref{defn:optimal-policy} clearly exists for well-defined problems, and its optimality can be proven in two logic progressions.
The first progression proves that if any policy $\pi^*$ satisfies \eqref{defn:optimal-policy}, then it also satisfies \eqref{eqn:value-equivalence}, which is done by backward induction.
As the initial step, \eqref{eqn:value-equivalence} trivially holds for $t=T+1$.
For the induction step, assume \eqref{eqn:value-equivalence} holds for $t=\tau+1$, and show it also holds for $t=\tau$. Recall \eqref{defn:Dt(h)} and proceed as follows
		\begin{eqnarray}
			D_{\tau-1}(h_{\tau};\beta) &=&  \min_{x_\tau\in X_\tau(h_\tau)} \max_{\zeta_\tau \in \mathcal{U}_\tau(\zeta_{1:\tau})} D_{\tau}(h_{\tau+1};\beta) \nonumber \\
			&=&  \max_{\zeta_\tau \in \mathcal{U}_\tau(\zeta_{1:\tau})} D_{\tau}(h_{\tau+1};\beta) \nonumber\\ 			
			&=&  \max_{\zeta_\tau \in \mathcal{U}_\tau(\zeta_{1:\tau})} D_{\tau}^{\pi^*}( h^{\pi^*}_{\tau+1}; \beta) \nonumber\\ 
			&=&  D_{\tau-1}^{\pi^*}( h_{\tau}; \beta), \nonumber
		\end{eqnarray}
where the second equality comes by \eqref{defn:optimal-policy}, the third by the inductive assumption, and the last comes by \eqref{defn:D[t](pi)}. Thus $\pi^*$ also satisfies \eqref{eqn:value-equivalence}.

The next progression proves that if a policy $\pi^*$ satisfies \eqref{eqn:value-equivalence}, then it is optimal under the principle of optimality, which becomes 
		\begin{equation} \label{eqn:value-dominance}
			D_{t-1}(h_{t};\beta) \leq D_{t-1}^\pi (h_{t};\beta), \forall \pi \in \Pi, \forall h_t\in H_t, t=1, \cdots, T+1,
		\end{equation}
	% for $t=1, \cdots, T+1$ via backward induction on $t$. 
after integrating \eqref{eqn:value-equivalence}. 
It is again by backward induction, whose initial step trivially holds for $t=T+1$.
For the induction step, assume (\ref{eqn:value-dominance}) holds for $t=\tau+1$ so as to show it also holds for $t=\tau$.
Recall (\ref{defn:Dt(h)}) with $t=\tau$ and apply the assumption by replacing $D_{\tau}(h_{\tau+1}; \beta)$ with $D_{\tau}^\pi (h_{\tau+1}; \beta)$:
\begin{eqnarray}
	D_{\tau-1}(h_{\tau};\beta) &\leq& 
		\min_{x_\tau\in X_\tau(h_\tau)} 
		\max_{\zeta_\tau \in \mathcal{U}_\tau(\zeta_{1:\tau})} 
		D_{\tau}^\pi (h_{\tau+1}; \beta) \nonumber \\
	&\leq& \max_{\zeta_\tau \in \mathcal{U}_\tau(\zeta_{1:\tau})} D_\tau^\pi (h^{\pi}_{\tau+1}; \beta) \nonumber\\ 
	&=&  D_{\tau-1}^\pi(h_{\tau}; \beta), \nonumber
\end{eqnarray}
where the second inequality comes by fixing $x_\tau = \pi_\tau(h_\tau)$, and the last line comes by (\ref{defn:D[t](pi)}).
Therefore (\ref{eqn:value-dominance}) holds by backward induction, and $\pi^*$ is optimal.
		 
%The reverse also holds. 
Once an optimal policy $\pi^*$ is known to satisfy \eqref{eqn:value-equivalence}, the principle of optimality requires that any optimal policy must satisfy \eqref{eqn:value-equivalence}, thus \eqref{defn:optimal-policy} automatically holds.
\end{proof}

\section{Proof of Theorem \ref{thm:1way-AR}}\label{app:1way-AR}
\begin{proof} 
By backward induction. For the initial step with $t=T$, it is easily verified. 
For the induction step, assume (\ref{eq:Dt}) holds in period $t+1\le T$ with 
\[ D_{t}(h_{t+1}; \beta) = \beta\max(\hat{p}_{t+1}, P_{T-t}(q_{t+1}))-R_{t+1},\]
and prove it also holds in period $t$.
For the minimization nested in \eqref{1way:D[t](h)}, let 
\begin{eqnarray}
	\bar{D}_t(h_{t}, p_t; \beta) &=& \min_{x_t \in X_t(h_t)}  D_t(h_{t+1}; \beta) \nonumber\\
	\label{defn:1way-barD[t]}
	&=& \min_{q_{t+1} \in [0, q_t]} \beta\max(\hat{p}_{t+1}, P_{n}(q_{t+1}))-R_{t+1},
\end{eqnarray}
with $n=T-t$, $q_{t+1} = q_t - x_t$, and $R_{t+1} = r_{t+1} + m q_{t+1}$.
% To find $\partial D_t(h_{t+1}; \beta)/\partial q_{t+1}$,
First note that
\[ P'_n( q ) = -{M-m\over\beta}\left(1-{q\over\beta n}\right)_+^{(n-1)}\leq 0, \]
which means $P_n(q)$ is monotone and there is
$P_n(q_{t+1}) \geq \hat{p}_{t+1}$  if $q_{t+1} \leq P^-_n(\hat{p}_{t+1})$ and $P_n(q_{t+1}) \leq \hat{p}_{t+1}$ otherwise. 
Therefore,
\begin{eqnarray}\label{eq:branchDt}
		D_t(h_{t+1}; \beta) &=& \left\{ \begin{array} {ll}
				\beta P_{n}(q_{t+1})-R_{t+1} & q_{t+1} \leq P^-_{n}(\hat{p}_{t+1})  \\
				\beta \hat{p}_{t+1} - R_{t+1}  &  q_{t+1} > P^-_{n}(\hat{p}_{t+1}) 
			\end{array} \right.\\
			%where $\bar{q}_t = g^-_t(\hat{p}_{t+1})$,  
			%\begin{equation} \label{eq:barqt}
			%\bar{q}_t = \beta t - \beta t
			%\left({\hat{p}_{t+1}-m \over M-m}\right)^{1/t}.
			%\end{equation}
			\label{eq:branch-pDt/pq}
		{\partial D_t(h_{t+1}; \beta) \over \partial q_{t+1}} &=& \left\{
			\begin{array}{ll}
				p_{t}-m +\beta P'_n(q_{t+1}) & q_{t+1} < P^-_{n}(\hat{p}_{t+1})  \\
				p_{t}-m  & q_{t+1} > P^-_{n}(\hat{p}_{t+1})
			\end{array} \right.
\end{eqnarray}
Note that with $q_{t+1} < P^-_{n}(\hat{p}_{t+1})$, there is 
		\(
		p_{t} \leq \hat{p}_{t+1} < P_n(q_{t+1}) 
		\leq - \beta P'_n(q_{t+1}) + m
		\), 
so $p_{t}-m+\beta P'_n(q_{t+1}) < 0$. 
And with $q_{t+1} > P^-_{n}(\hat{p}_{t+1})$, there is $p_{t}-m \geq 0$. 
It is clear that (\ref{eqn:1way-q*-solved}) is an optimal solution to (\ref{defn:1way-barD[t]}), which from (\ref{eq:branchDt}) gives 
\begin{equation}  \label{eqn:1way-barDt-solved}
			\bar{D}_t( h_t, p_t ; \beta) = \beta P_{n}(q_{t+1}^*)-(r_{t+1}+mq_{t+1}^*).
\end{equation}
Let $\bar{p}_t=\max(\hat{p}_{t}, P_n(q_t)) \in [m,M]$, and from (\ref{1way:D[t](h)}) there is 
\begin{eqnarray} 
	D_{t-1}( h_t; \beta) &=& \max_{p_{t}\in[m,M]} \bar{D}_t( h_t, p_t ; \beta) \nonumber\\
	\label{max:split}
	&=& \max\left(
	\begin{array}{c}
		\max_{p_{t}\in[m,\bar{p}_t]} \bar{D}_t( h_t, p_t ; \beta)  \\
		\max_{p_{t}\in[\bar{p}_t,M]} \bar{D}_t( h_t, p_t ; \beta)
	\end{array} 
	\right) 
\end{eqnarray}
		
For the branch with $p_t \in [m,\bar{p}_t]$ in (\ref{max:split}), consider two cases: (i) $\bar{p}_t = \hat{p}_{t} \geq P_n(q_t)$ and (ii) $\bar{p}_t = P_n(q_t) > \hat{p}_{t}$. 
In case (i) there is $\hat{p}_{t+1} = \max(\hat{p}_t, p_t) = \hat{p}_t \geq P_n(q_{t})$, therefore $P^-_n(\hat{p}_{t+1}) \leq q_t$ and (\ref{eqn:1way-q*-solved}) simplifies to $q^*_{t+1}  = P^-_n(\hat{p}_{t+1})$, thus $P_n(q^*_{t+1}) = \hat{p}_{t+1} = \bar{p}_t$.
In case (ii) there is $\hat{p}_{t+1} \leq P_n(q_{t})$, therefore $P^-_n(\hat{p}_{t+1}) \geq q_t$ and  (\ref{eqn:1way-q*-solved}) simplifies to $q^*_{t+1} = q_t$, thus $P_n(q^*_{t+1}) = P_n(q_t) = \bar{p}_t$. 
So there is $P_n(q^*_{t+1}) = \bar{p}_t$ in both cases, and (\ref{eqn:1way-barDt-solved}) becomes $\bar{D}_t( h_t, p_t ; \beta) = \beta \bar{p}_t - (r_{t+1}+m q_{t+1}^*) = \beta \bar{p}_t - r_t - p_t x^*_t - m q_{t+1}^*$, which is linear in $p_{t}$ with a slope of $- x^*_t \leq 0$ as $x^*_t = q_{t}-q_{t+1}^*\geq 0$. Thus $p^*_{t}=m$ is a maximizer, which gives
\(
\max_{p_{t}\in[m,\bar{p}_t]} \bar{D}_t( h_t, p_t ; \beta) = \beta \bar{p}_t - r_t - m q_{t} = \beta \bar{p}_t - R_t.
\)
		
For the branch with $p_t \in [\bar{p}_t, M]$ in (\ref{max:split}), as $p_t \geq \bar{p}_t \geq \hat{p}_t$, there is $\hat{p}_{t+1} = p_{t} \geq \bar{p}_t \geq P_n(q_{t})$, thus $P^-_n(\hat{p}_{t+1}) \leq q_t$ and (\ref{eqn:1way-q*-solved}) simplifies to $q^*_{t+1} = P^-_n(\hat{p}_{t+1})$. 
Therefore $P_n(q^*_{t+1}) = \hat{p}_{t+1} = p_t$, and (\ref{eqn:1way-barDt-solved}) simplifies to
$\bar{D}_t( h_t, p_t ; \beta) = \beta p_t - r_t - p_t x^*_t - m q_{t+1}^* = \beta p_t - r_t - p_t (q_t - q_{t+1}^*) - m q_{t+1}^* = (\beta - q_t + q_{t+1}^*)  p_{t}  - m q_{t+1}^*- r_t = (\beta - q_t + q_{t+1}^*)  P_n(q^*_{t+1}) - m q_{t+1}^*  - r_t = d(P^-_n(p_t))$, where
\(
d(z) = (\beta - q_t + z)  P_n(z) - m z  - r_t, z\in[0, 1], %z\in[0, q_t]
\) 
with a derivative $d'(z) = (\beta - q_t + z)P'_n(z) + P_n(z) - m$. 
Note that $P_n(z) - m = - (\beta-z/n) P'_n(z)$, thus $d'(z) = (\beta - q_t + z)P'_n(z) - (\beta-z/n) P'_n(z) = (z+z/n - q_t) P'_n(z)$.
As $P'_n(z) \leq 0$, there is $d'(z) \geq 0$ when $z+z/n - q_t \leq 0$, and $d'(z) \leq 0$ when $z+z/n - q_t \geq 0$, hence $z^* = n q_t / (n+1) < q_t$ solves $\max_{z\in [0,1]} d(z)$, which gives  
\[
		d(z^*)  
		% &=&(\beta - q_t + z^*) P_n(z^*) - m z^* - r_t\\
		% &=&  (\beta - q_t + z^*) (M-m)\left(1-{z^*\over
		% \beta n}\right)_+^{n} + \beta m - m q_t - r_t\\
		% &=& \beta  (1 - {q_t\over \beta(n+1)}) (M-m) 
		% \left(1-{q_t\over \beta (n+1)}\right)_+^{n} + 
		% \beta m - R_t\\
		% &=&\beta (M-m) \left(1-{q_t\over \beta 
		% (n+1)}\right)_+^{(n+1)} + \beta m - R_t\\
		= \beta P_{n+1}(q_t) - R_t, ~
		P_n(z^*) \ge  P_{n+1}(q_t).
\]
Consider two cases with $\bar{D}_t( h_t, p_t ; \beta) =
	d(P^-_n(p_t))$ for $p_{t}\in[\bar{p}_t,M]$.
Case (i) $P_n(z^*) \geq \bar{p}_t$. 
As $P^-_n(M)=0 \leq z^* \leq P^-_n(\bar{p}_t)$,
there is $\max_{p_{t}\in[\bar{p}_t,M]} \bar{D}_t( h_t, p_t ; \beta) = d(z^*)$. 
Thus, according to (\ref{max:split}) there is 
\begin{equation}\label{eqn:D[t-1]-simplified}
D_{t-1}( h_t; \beta) = \max (\beta \bar{p}_t - R_t, d(z^*)).
\end{equation} 
Case (ii) $P_n(z^*) < \bar{p}_t$. 
As $q_t \geq z^* \geq P^-_n(\bar{p}_t)$,
there is 
\begin{eqnarray*}
\max_{p_{t}\in[\bar{p}_t,M]} \bar{D}_t( h_t, p_t ; \beta) &=& \max_{p_{t}\in[\bar{p}_t,M]} d(P^-_n(p_t)) \\
&=& \max_{z\in [0, P^-_n(\bar{p}_t)]} d(z)\\
&\leq& \max_{z\in [0, q_t]} d(z) = d(z^*).	
\end{eqnarray*}
As $P_n(z^*) \geq P_{n+1}(q_t)$,
there is $\bar{p}_t \geq P_{n+1}(q_t)$.
So $d(z^*) = \beta P_{n+1}(q_t) - R_t \leq \beta \bar{p}_t - R_t$, and by (\ref{max:split}) there is $D_{t-1}( h_t; \beta) = \beta \bar{p}_t - R_t$, and \eqref{eqn:D[t-1]-simplified} remains valid.
Therefore, in both cases proceed from \eqref{eqn:D[t-1]-simplified} and take note of $\bar{p}_t = \max(\hat{p}_{t}, P_n(q_{t}))$ and $P_n(q_{t}) \leq P_{n+1}(q_{t})$:
\begin{eqnarray*}
		D_{t-1}( h_t; \beta) & = & \max (\beta \bar{p}_t - R_t, d(z^*))\\
		&=& \max (\beta \bar{p}_t - R_t, \beta P_{n+1}(q_t) - R_t) \\
		&=& \beta\max ( \bar{p}_t , P_{n+1}(q_t)) - R_t\\
		&=& \beta \max(\hat{p}_{t}, P_n(q_{t}), P_{n+1}(q_{t})) -R_t\\
		&=&\beta \max(\hat{p}_{t},  P_{n+1}(q_{t})) - R_t
\end{eqnarray*}
As $n=T-t$, clearly (\ref{eq:Dt}) also holds for $t$. 
\end{proof} 

%\bibliographystyle{elsarticle-harv} 

\bibliographystyle{apalike} 
\bibliography{cas-refs}

%\input{cas-refs.bbl}





% \end{thebibliography}
\end{document}